\documentclass[12pt,a4paper]{article}
\pdfoutput=1
\usepackage{jheppub}
 \usepackage{amsmath, amssymb, latexsym, amsthm, amsfonts, mathrsfs, amscd, mathtools}
 \usepackage{xcolor}
 \usepackage{caption}
 \usepackage{subcaption}
 \usepackage{slashed}  
\usepackage{float}
\usepackage{graphicx}
\graphicspath{ {project2 tex/} }

\def\CE{{\mathcal{E}}}
\def\CM{{\mathcal{M}}}

\def\CR{{\mathcal{R}}}

\def\b{\bullet}

\newcommand{\be} {\begin{equation}}
\newcommand{\ee} {\end{equation}}
\newcommand{\bea} {\begin{eqnarray}}
\newcommand{\eea} {\end{eqnarray}}
\newcommand{\ba} {\begin{array}}
\newcommand{\ea} {\end{array}}
\newcommand{\nn} {\nonumber}

\newcommand{\ul}{\underline}

\theoremstyle{plain}

\theoremstyle{definition}

\newcommand{\sgn} { {\rm sgn}}

 \title{Appell Functions for General Lattices}
 \author{Aradhita Chattopadhyaya$^a$, Jan Manschot$^{b,c}$}
\affiliation{ $^a$
Chennai Mathematical Institute,\\
H1 SIPCOT IT Park\\
Siruseri, Kelambakkam\\
Tamil Nadu 603103, India\\
$^b$School of Mathematics, Trinity College, Dublin 2,
   Ireland\\
   $^c$Hamilton Mathematical Institute, Trinity College, Dublin 2,
   Ireland\\}
\emailAdd{aradhitac@cmi.ac.in, jan.manschot@tcd.ie}
  
\abstract{We study Appell functions associated to an arbitrary positive definite lattice $\Lambda$ and a
  choice of $M\leq {\rm
  dim}(\Lambda)$ linearly independent vectors $d_r\in \Lambda$, $r=1,\dots,M$. These
functions are instances of multi-variable
quasi-elliptic functions, and specific examples have appeared at various places in
mathematics and theoretical physics. For example, if $\Lambda$ is chosen to be one-dimensional, these functions reduce to the classical Appell
function, which is a prominent example in the theory of mock modular
forms. The Appell functions introduced here are examples of depth $M$ mock modular forms. We derive a structural formula for their modular completion. Motivated by partition functions in theoretical physics,
we discuss the case where $\Lambda$ is the $A_N$ root lattice in detail. \\
\\
}    
\begin{document} 
\maketitle       
\section{Introduction}  
\label{sec:intro} 
The Appell function is a classic two-variable quasi-elliptic function, whose
history goes back to the 19'th century \cite{Appell:1886,
  Lerch:1892}.\footnote{This elliptic Appell function and the
  functions studied in this paper are not to be confused with the Appell
  hypergeometric functions $F_j$.} 
In their long history, specializations of Appell functions were
identified as examples of Ramanujan's mock theta functions \cite{ramanujan_lost_notebook, Watson:1936}, and now play a central role
in the theory of mock modular forms 
\cite{ZwegersThesis, MR2605321}. These functions have found applications across many disciplines in
mathematics and theoretical physics, including conformal field theory \cite{Eguchi:1988af, Kac:2001,
   semikhatov2005, Cheng:2012rca}, algebraic geometry
\cite{polishchuk2001}, supersymmetric black holes \cite{Dabholkar:2012nd} and
topological quantum field theory \cite{Bringmann:2010sd, Manschot:2014cca,
  Korpas:2019cwg, Manschot:2021qqe}. This has motivated in part the study of
variants of the original Appell functions such as those with
multiple elliptic variables in Ref. \cite{Zwegers:2019}. The present
paper will study a general family of Appell functions based on a
positive definite lattice $\Lambda$ and $M\leq \dim(\Lambda)$ vectors $d_r\in
\Lambda$. These are examples of mock Jacobi forms with depth $M$.  
Higher depth mock modular forms have found applications in mathematical physics
\cite{Manschot:2014cca, Manschot:2017xcr, KumarGupta:2018rac, Alexandrov:2019rth, Alexandrov:2018lgp,
  Chattopadhyaya:2021rdi, Alexandrov:2024jnu, Alexandrov:2024wla}. We will further elaborate on this connection below. 
 
To state the Appell functions studied in this paper, let $\Lambda$ be
a positive definite $N$-dimensional lattice, with quadratic form $Q: \Lambda\to \mathbb{Z}$
and bilinear form $B: \Lambda \times \Lambda \to \mathbb{Z}$. Let
$\{d_r\}$, $r=1,\dots, M\leq N$,  be a set of $M$ linearly
independent vectors $d_r\in \Lambda$. With $\mu\in 
\Lambda \otimes \mathbb{R}$ and $\mathbb{H}$ the upper-half plane, this article studies the 
Appell function,
\be 
\Phi_{\mu}:\mathbb{H}\times (\Lambda \otimes \mathbb{C})^2 \to \mathbb{C},
\ee
defined as \cite{Manschot:2014cca}
\begin{eqnarray}\label{phigen1} 
\Phi_{\mu}(\tau,u,v,\{ d_r \}) &=& \sum_{k\in{\Lambda+{\mu}}}
                                       \frac{q^{Q(k)/2}e^{2\pi
                                       iB({v},k)}}{\prod_{r=1}^M(1-e^{2\pi
                                       i B(d_r,u)} q^{B(d_r,
                                       k)})},    
\end{eqnarray}
where $q=e^{2\pi i \tau}$. For $\Lambda\simeq \mathbb{Z}$ and $M=1$, $\Phi_{\mu}$ reduces to a
variant of the classical Appell function. We set out to study the more general Appell
function $\Phi_{\mu}$ using techniques for
indefinite theta series as employed earlier \cite{Alexandrov:2016enp,
  Nazaroglu:2016lmr, Manschot:2017xcr}. This makes it clear that
$\Phi_{\mu}$ is an example of a mock
modular form or mock Jacobi form of depth $M$. That is to say,
$\Phi_{\mu}$ does not transform in the standard way under modular
transformations. However for a specific non-holomorphic function
$\CR_\mu$, the non-holomorphic completion,
\be
\widehat
\Phi_{\mu}(\tau,u,v,\{d_r\})=\Phi_{\mu}(\tau,u,v,\{d_r\})+\CR_\mu(\tau,\bar
\tau,u,\bar u, v, \bar v,\{d_r\})
\ee 
does transform as modular or Jacobi form. For a mock modular form of
depth $M$, $\CR_\mu$ involves involves an $M$-dimensional iterated
integral. We provide an explicit structural formula
(\ref{eq:StructForm}) for the modular completion $\widehat
\Phi_{\mu}$. It establishes that the completion can be compactly written in
terms of Appell functions for the same lattice $\Lambda$, but with
depth $M'<M$ multiplied by non-holomorphic iterated integrals $M_L$ \cite{Alexandrov:2016enp}. The
latter integrals can be seen as a higher-dimensional generalization of
the complementary error function.

To this end, we relate the $N$-dimensional lattice $\Lambda$ and set
$\{d_r\}$ to an $(M+N)$-dimensional lattice $\ul \Lambda$ of signature
$(N,M)$, and apply techniques from indefinite theta series
\cite{Vigneras:1977, ZwegersThesis, Alexandrov:2016enp,
  Nazaroglu:2016lmr}. The $N$-dimensional elliptic variable $u$ and
$M$-dimensional elliptic variable $v$ combine to an
$(M+N)$-dimensional variable $\ul z \in \ul \Lambda \otimes \mathbb{C}$. An attractive feature of the Appell functions
compared to generic indefinite theta functions is that the Appell
functions depend on $M$ vectors while the indefinite theta series
would depend on $2M$. This reduces the complexity of their
non-holomorphic terms significantly. 
   
One of the motivations of the authors to study these functions, is
their appearance in topologically twisted, $\mathcal{N}=4$ supersymmetric Yang-Mills
theory and related algebraic geometry \cite{Vafa:1994tf,
  Yoshioka94, Yoshioka:k95, Bringmann:2010sd, Manschot:2011ym, 
  Manschot:2011dj, Bringmann_2016, toda2017,    
  Bringmann:2018cov, Alexandrov:2020bwg, Alexandrov:2020dyy}. The
physical partition functions of topologically twisted gauge theory
with gauge group $SU(N)$ exhibit a holomorphic anomaly, which is of much
interest in physics \cite{Vafa:1994tf, Minahan:1998vr, Dabholkar:2020fde,
   Alexandrov:2019rth}. The proposal that the holomorphic anomaly
 involves the partition functions for  groups $SU(N')$ with $N'<N$
 \cite[Eq. (4.7)]{Minahan:1998vr} is confirmed in many cases. The generalized Appell function $\Phi_{\mu}$ (\ref{phigen1}) arises as
 building blocks of the generating functions of Poincar\'e
polynomials of moduli spaces of sheaves, which are derived using
Harder-Narasimhan filtrations and algebraic-geometric invariants for
stacks \cite{Yoshioka94, joyce2008, Manschot:2014cca}.  The main building block of the partition functions derived in
\cite{Manschot:2014cca} for gauge group $SU(N)$ are $\Psi_{(r_1,\dots,
  r_{\ell})(a,b)}$ with $\sum_j r_j=N$, defined
as:
\begin{eqnarray}\label{psi} 
\Psi_{(r_1,\dots ,r_{\ell}),(a,b)}(\tau,z) &=&
                                             \sum_{\begin{smallmatrix}\sum_{s=1}^{\ell}b_sr_s=b\\b_j\in\mathbb{Z}\end{smallmatrix}} \frac{w^{\sum_{j<i}r_ir_j(b_i-b_j)+\sum_{i}2(r_i+r_{i-1})\left\{\frac{a}{r}\sum_{k=i}^{\ell} r_k\right\}}}{\prod_{i=2}^{\ell}(1-w^{2(r_i+r_{i-1})}q^{b_{i-1}-b_i})}\\ \nn
&& \times q^{\sum_{i=1}^{\ell}\frac{r_i(r-r_i)}{2r}b_i^2-\frac{1}{r}\sum_{i<j}b_ib_jr_ir_j+\sum_{i=2}^{\ell}(b_{i-1}-b_i)\left\{\frac{a}{r}\sum_{k=i}^{\ell}r_k\right\}},
\end{eqnarray}
where $w=e^{2\pi i z}$, and $\{\,\}:\mathbb{R}\to [0,1)$ is the fractional part,
\be
\{x\}=x-\lfloor x\rfloor,\qquad \{-x\}=-x+\lceil x \rceil,
\ee 
$r_i\in \mathbb{N}^*$, $r=\sum_{i=1}^\ell r_i$, $a,b\in 
\mathbb{Z}$. If $r_i=1$, for all $i=1,\dots,\ell$, the associated
quadratic form is $A_{\ell-1}$. These can indeed be identified as
examples of $\Phi_{\mu}$ or the slightly more general $\Phi_{\mu,\nu}$
(\ref{phigen2}). The analysis of this paper will make it more
straightforward to determine and write the completion of functions
such as $\Psi_{(r_1,\dots
  r_{\ell}),(a,b)}$.  
  
Relatedly, the Appell functions may find applications
for partition functions of supersymmetric black holes, which are known
to involve mock modular forms \cite{Manschot:2007ha, Manschot:2009ia,
  Dabholkar:2012nd, Alexandrov:2016tnf, Alexandrov:2018lgp, Alexandrov:2023ltz}. In
particular, it could aid the determination of the holomorphic part of a
partition function given its non-holomorphic part \cite{Alexandrov:2017qhn, Alexandrov:2024wla}.  
    
 The outline of this paper is as follows. Section \ref{sec:AppellFs}
 introduces the Appell function $\Phi_\mu$ for a general lattice $\Lambda$ and
 depth $M\leq {\rm dim}(\Lambda)$. This section also develops the
 connection with an indefinite theta series for a lattice $\ul
 \Lambda$. Section \ref{sec:ModCom} determines the modular completion
 $\widehat \Phi_\mu$ and derives the structural formula
 \eqref{eq:StructForm}. Section \ref{sec:AppellBPS} applies the
 general formulas to Appell functions appearing in the context of BPS
 indices.

\subsection*{Acknowledgements}
We thank Sergey Alexandrov and Caner Nazaroglu for useful discussions
and correspondence. The work of A.C. is supported by Ramanujan
Fellowship RJF/2023/000070 from the Anusandhan National Research Foundation, India.

 \section{Appell Functions}
\label{sec:AppellFs}
 We introduce in this section the general family of Appell functions
for an arbitrary positive definite lattice $\Lambda$. 
In Subsection \ref{sec:Anlattice}, we will specialize to the case where $\Lambda$ is
the root lattice of $SU(N+1)$. 

Throughout, we let $\tau\in \mathbb{H}$, $y={\rm Im}(\tau)$ and $q=e^{2\pi i \tau}$.

\subsection{Preliminaries on Lattices}
As in Section \ref{sec:intro}, we consider an $N$-dimensional positive
definite, integral lattice $\Lambda$, with
bilinear form $B_\Lambda=B$ and quadratic
form $Q_\Lambda=Q$ related through $Q(k)=B(k,k)$. The dual lattice of $\Lambda$ is denoted by
$\Lambda^*$. For a positive definite
lattice $\Lambda$, we let $-\Lambda$ be the lattice with negative
definite bilinear form $B_{-\Lambda}=-B_{\Lambda}$.

In the following, we will often consider a set of $M$ linearly
independent vectors
$\{d_r\}$, $d_r\in \Lambda$ with $r=1,\dots, M$. We define the sublattice $\Lambda(\{d_r\})\subseteq
\Lambda$ as the $M$-dimensional sublattice of $\Lambda$ generated by
$\{d_r\}$. If there is no confusion on the set $\{d_r\}$, we will
sometimes abbreviate $\Lambda(\{d_r\})$ to
$\Lambda_d$.

\subsection{Definition and First Properties} 

We introduce here a slight variation of the Appell function
$\Phi_{\mu,\nu}$ compared to $\Phi_{\mu}$ \eqref{phigen1}. 
As above, we choose a set $\{d_r\}$ of $M$ linearly independent
vectors $d_r\in \Lambda$ spanning $\Lambda_d$. Furthermore, let
$u\in\Lambda_d\otimes\mathbb{C}$, $v\in \Lambda \otimes \mathbb{C}$,
$\mu\in \Lambda\otimes \mathbb{R}$, and $\nu\in \Lambda_d\otimes
\mathbb{R}$. We will often consider $\mu\in \Lambda^*$ and $\nu\in \Lambda_d^*$.

We then define the Appell function $\Phi_{\mu,\nu}$ as\footnote{This
  definition is a variation on the definitions in
  \cite[Eq. (0.13)]{kac2014} and
  \cite[Eq. (5.2)]{Manschot:2014cca}. An important difference with the
  definition in Ref. \cite{kac2014} is that the set $\{d_r\}$ in
  Eq. (\ref{phigen2}) is not required to consist of pairwise
  orthogonal vectors.}
\be
\begin{split}
\label{phigen2} 
\Phi_{\mu,\nu}(\tau,u,v,\{ d_r \}) &= e^{2\pi i B(\nu,u-v)}q^{-Q(\nu)/2}\sum_{k\in{\Lambda+{\mu}}}
                                       \frac{q^{Q(k)/2}e^{2\pi
                                       iB({v},k)}}{\prod_{r=1}^M(1-e^{2\pi
                                       i B(d_r,u)} q^{B(d_r,
                                           k-\nu)})} \\
  &= e^{2\pi i B(\nu,u)}\sum_{k\in{\Lambda+{\mu-\nu}}}
                                       \frac{q^{Q(k)/2+B(\nu,k)}e^{2\pi
                                       iB({v},k)}}{\prod_{r=1}^M(1-e^{2\pi
                                       i B(d_r,u)} q^{B(d_r,
                                       k)})}.
\end{split}
\ee   
This function $\Phi_{\mu,\nu}$ can be expressed in terms of
$\Phi_{\mu}$, since $\Phi_{\mu,0}=\Phi_{\mu}$  and
\be
\label{eq:PhimunuPhimu0}
\Phi_{\mu,\nu}(\tau,u,v,\{ d_r \})=e^{2\pi i
  B(\nu,u)}\Phi_{\mu-\nu,0}(\tau,u,v+\nu\tau, \{ d_r \}).
\ee

In analogy with the original Appell function, $\Phi_{\mu,\nu}$ has two
elliptic arguments $u$ and $v$.  For fixed argument $\tau\in \mathbb{H}$,
$\Phi_{\mu,\nu}(\tau,u,v,\{d_r\})$ is holomorphic in $v$,
while meromorphic in $u$ with poles at 
\be
B(d_r,u+k\tau)\in \mathbb{Z}, \quad {\rm for\,\, all} \quad k\in \Lambda+\mu-\nu.
\ee

It is straightforward to check that the following (quasi)-periodicity properties hold:
\begin{enumerate}
\item For the shift $\mu \mapsto \mu+\ell$ with $\ell\in \Lambda$:
\be
\Phi_{\mu+\ell,\nu}(\tau,u,v,\{ d_r \})=\Phi_{\mu,\nu}(\tau,u,v,\{ d_r \}).
\ee
For the shift $\nu\mapsto \nu +\ell_d$ with $\ell_d\in \Lambda_d$:
\be
\label{eq:nuqperiod}
\Phi_{\mu,\nu+\ell_d}(\tau,u,v,\{ d_r \})=e^{2\pi i B(\ell_d,u)}\Phi_{\mu,\nu}(\tau,u,v+\ell_d\tau,\{ d_r \}).
\ee 
For the simultaneous shift of $\mu$ and $\nu$ by $m_d\in \Lambda_d^*$:
\begin{eqnarray}
\Phi_{\mu+m_d,\nu+m_d}(\tau,u,v,\{ d_r \}) &=& e^{2\pi i B(m_d,u)}\Phi_{\mu,\nu}(\tau,u,v+m_d\tau,\{ d_r \}).
\end{eqnarray}
Furthermore for $m=m^{||}+m^{\perp}\in \Lambda^*$ with $m^{||}$
and $m^{\perp}$ the components parallel and orthogonal to $\Lambda_d$, we have
\be
\label{eq:Phimum}
\begin{split}
  \Phi_{\mu+m,\nu+m^{||}}(\tau,u,v,\{d_r\})&=e^{2\pi i B(m^{||},u)+2\pi i
  B(v,m^\perp)}q^{Q(m^\perp)/2}
\\
&\quad \times \Phi_{\mu,\nu}(\tau,u,v+m\tau,\{d_r\}). 
\end{split}
\ee

\item For the inversion $(\mu,\nu)\to -(\mu,\nu)$:
  \be
\label{Phi-mu}
  \begin{split}
    &\Phi_{-\mu,-\nu}(\tau,u,v,\{ d_r \})= \\
    &\qquad (-1)^M e^{-2\pi i B(\sum_r
    d_r,u)} \Phi_{\mu,\nu}(\tau,-u,-v+\tau \sum_r d_r,\{
  d_r \}).
  \end{split}
  \ee
\item For the shift of $u$ by $m_d\in\Lambda_{d}^*$:
\begin{eqnarray}\label{rel1}
\Phi_{\mu,\nu}(\tau,u+m_d,v,\{ d_r \})=e^{2\pi i B(m_d,\nu)}\Phi_{\mu,\nu}(\tau,u,v,\{ d_r \}),
\end{eqnarray}
since $B(d_r,m_d)\in\mathbb{Z}$.

For a shift of $v$ by $m\in\Lambda^*$:
\begin{eqnarray}
\Phi_{\mu,\nu}(\tau,u,v+m,\{ d_r \})=e^{2\pi i B(m,\mu-\nu)} \Phi_{\mu,\nu}(\tau,u,v,\{ d_r \}), 
\end{eqnarray}
since $B(k,m)= B(m,\mu-\nu) \mod \mathbb{Z}$.

For the simultaneous shift of $u$ and $v$ by $\ell_d\tau$ with $\ell_d\in\Lambda_{d}\subseteq \Lambda$:
\be\Phi_{\mu,\nu}(\tau,u+\ell_d\tau,v+\ell_d\tau,\{ d_r
\})=q^{-Q(\ell_d)/2}e^{-2\pi i
  B(v,\ell_d)}\Phi_{\mu,\nu}(\tau,u,v,\{ d_r \}).
\ee

\end{enumerate}

Before moving to the next section, we consider the action of a
matrix $G\in \mathrm{SL}(N,\mathbb{Z})$ on $k\in \Lambda$, which leaves invariant the
bilinear and quadratic form
\be
\label{eq:TrafoInv}
B(Gk,Gk')=B(k,k'),\qquad  Q(Gk)=Q(k).
\ee
While this transformation leaves invariant the lattice $\Lambda$, it
will in general transform the sublattice $\Lambda_d$ to a
different sublattice $\Lambda_{d'}$ generated by $d_r'=Gd_r$. This transformation acts on
$\Phi_{\mu,\nu}$ as follows
\begin{eqnarray}
\label{eq:PhiGInv}
  \Phi_{\mu,\nu}(\tau,u,v,\{ d_r \})&=& e^{2\pi i B(Gu,G\nu)} \sum_{k\in{\Lambda+{\mu-\nu}}}
\frac{q^{Q(G k)/2+B(G\nu,Gk)}e^{2\pi iB(G{v},Gk)}}{\prod_{r=1}^M(1-
    e^{2\pi i B(Gd_r,G u)} q^{B(Gd_r,G k)} )}\\ \nn
&=& e^{2\pi i B(Gu,G\nu)} \sum_{k\in{\Lambda+{G(\mu-\nu)}}}
\frac{q^{Q(k)/2+B(G\nu, k)}e^{2\pi iB(G{v},k)}}{\prod_{s=1}^M(1-e^{2\pi i B(d_r',G u)} q^{B(d_r',k)} )}\\ \nn
&=& \Phi_{G\mu,G\nu}(\tau, Gu, Gv,\{ d'_r \}).
\end{eqnarray}

If the set $\{d_r\}$ is the empty set $\varnothing$, $\Phi_{\mu,\nu}$
is independent of $u$ and $\nu$. It is simply a theta series for the
lattice $\Lambda$. Although not widely used in this paper, it is convenient to also
introduce the normalized Appell function,
\be
M_{\mu,\nu}(\tau,u,v,\{d_r\})=\frac{\Phi_{\mu,\nu}(\tau,u,v,\{ d_r \})}{\Phi_{\mu}(\tau,-,v,\varnothing)},
\ee
where $-$ indicates that the
function is independent of $u$.
This is the higher-dimensional analogue of
the Lerch sum $\mu(\tau,u,v)$ of Ref. \cite{ZwegersThesis}.
It is intriguing that the coefficients of these functions exhibit
moderate growth in many examples \cite{Chattopadhyaya:2023aua}. That is to say the function is weakly holomorphic as
function of $\tau$. 
  
\subsection{Appell Functions as Indefinite Theta Series}
\label{sec:AppellIndef}
  
In this subsection, we will relate the Appell function (\ref{phigen1}) to an indefinite theta series. This will be important
in Section \ref{sec:ModCom} to determine the modular completion $\widehat
\Phi_\mu$ of $\Phi_\mu$. We start by recalling
the definition of an indefinite theta series.

\subsubsection*{Indefinite Theta Series}
An indefinite theta series is a holomorphic
$q$-series obtained from a sum over an indefinite lattice
\cite{Gottsche:1996aoa, ZwegersThesis, Alexandrov:2016enp,  Nazaroglu:2016lmr, kudla2016theta}. 
Let $\Gamma$ be an $(M+N)$-dimensional
indefinite lattice of signature $(N,M)$ with quadratic form $Q_\Gamma=Q$ and
bilinear form $B_\Gamma=B$. 
We define the indefinite
theta series $\Theta_{\Gamma, \mu}:\mathbb{H}\times
(\Gamma \otimes \mathbb{C})\to \mathbb{C}$ as,   
\be
\label{eq:DefIndefTheta}
\Theta_{\Gamma,\mu}(\tau, z,\{C_r,C_r'\})=\sum_{k\in \Gamma+\mu}
K(\{C_r,C_r'\}, k+  a)\,q^{Q(k)/2}e^{2\pi i  B( k, z)},
\ee
where $ a={\rm Im}(z)/y$, and the support of the kernel $K$ is
such that the sum over $ k\in \Gamma+ \mu$
is convergent. We will consider $K$ of the form 
\be
\label{eq:KIndefTheta}
K(\{C_r,C_r'\},x)=2^{-M} \prod_{r=1}^{M} \left( \mathrm{sgn}(
  B(x,C_r)) + \mathrm{sgn}(B(x,C_r'))\right),
\ee
with $x\in \Lambda\otimes \mathbb{R}$. 
Convergence puts non-trivial constraints on the set $\{C_r,C_r'\}$
\cite{ZwegersThesis, Alexandrov:2016enp, Nazaroglu:2016lmr, kudla2016theta,
  2017arXiv170802969F, FunkeNotes:2018}. 

\subsubsection*{Appell Functions}
To relate $\Phi_\mu$ (\ref{phigen2}) to an indefinite theta series
$\Theta_{\Gamma,\mu}$ (\ref{eq:DefIndefTheta}), we
expand the denominator using a
geometric series expansion as
\begin{eqnarray}
\label{PhiSigns} 
\Phi_{\mu,\nu}(\tau, u, v,\{d_r\}) &=& e^{2\pi i B(\nu,u)}
                                 \sum_{x_r\in\mathbb{Z}}\sum_{k\in
                                 \Lambda +\mu-\nu}
                                 q^{\frac{1}{2}{Q}(k)+B(\nu,k)} e^{2\pi i B(v,k)}\\
&& \times 2^{-M}\prod_{r=1}^M ({\rm sgn}(x_{r}+\epsilon)+{\rm
   sgn}(B(d_r,k+a)))\,e^{2\pi i B(x_rd_r,u)}
   q^{{B}(x_rd_r,k)}\nn\\
&=& \sum_{x_r\in\mathbb{Z}}\sum_{k\in
                                 \Lambda +\mu-\nu}
    q^{\frac{1}{2}{Q}(k)+B(\nu+x_rd_r,k)}
    e^{2\pi i B(u,\nu+x_rd_r)+ 2\pi i B(v,k)}\nn \\
&&\times 2^{-M}\prod_{r=1}^M ({\rm sgn}(x_{r}+\epsilon)+{\rm
   sgn}(B(d_r,k+a))),\nn
\end{eqnarray}
with $0<\epsilon\ll 1$ and $a={\rm Im}(u)/y\in \Lambda \otimes
\mathbb{R}$. This expression demonstrates that $\nu+\sum_r
x_rd_r\equiv \sum_r (\nu_r+x_r)d_r$ lies
naturally in $\Lambda_d^*$, and the expression resembles the form of the indefinite
theta series (\ref{eq:DefIndefTheta}). To make the correspondence more
precise, we need to identify:
\begin{enumerate}
\item the indefinite lattice $\Gamma$ for the indefinite theta series,
\item the elliptic variable $z$ for the indefinite theta series,
\item the vectors $C_r$ and $C_r'$,
\item the kernel $K$. 
\end{enumerate}
 
\subsubsection*{The Lattice}
We let the lattice $\Lambda$ be spanned by the set of vectors $\alpha_j$, $j=1,\dots N$.
To determine the lattice $\Gamma$ for the indefinite theta series, we extend the lattice vector $k=\sum_j
k_j\alpha_j\in \Lambda$ with the $x_r$ to form an 
$(M+N)$-dimensional vector ${\underline k}\in \Gamma$. The following discussion will demonstrate that
$\Gamma$ equals $\Lambda\oplus
(-\Lambda_d)\subseteq \Lambda\oplus
(-\Lambda)$. We will denote the lattice $\Lambda\oplus
(-\Lambda_d)$ by $\ul \Lambda$ and distinguish elements in $\ul
\Lambda$ also with an underline, for example $\ul k$, $\ul \mu$ and
$\ul z$. 

The natural basis elements of $\ul \Lambda$ are $\ul \alpha_j$, $j=1,\dots,
M+N$ with $\ul \alpha_i=(\alpha_i,0)\in \underline
\Lambda\subset \Lambda\oplus 
(-\Lambda)$, $i=1,\dots, N$ together with $\ul \alpha_{N+r}=(0,d_r)\in \underline
\Lambda\subset \Lambda\oplus
(-\Lambda)$, $r=1,\dots,M$. The lattice $\ul \Lambda$ comes with a
quadratic form $\ul B: \ul \Lambda\times \ul \Lambda \to \mathbb{Z}$, which
evaluates on the basis elements $\ul \alpha_i$ as
\be
\begin{split} 
&\ul{B}(\ul \alpha_i,\ul \alpha_j)=B(\alpha_i,\alpha_j), \qquad \qquad \quad
i,j\in \{1,\dots,N\}, \\
&\ul{B}(\ul \alpha_{N+r},\ul \alpha_{N+s})=-B(d_r,d_s),\qquad
r,s\in \{1,\dots,M\},
\end{split} 
\ee
and else 0. Another useful basis is $\{\ul \alpha_j'\}$ with $\ul
\alpha_i'=\ul \alpha_i$ for $i=1,\dots,N$ and null vectors $\ul
\alpha_{N+r}'=\ul \gamma_r$, $r=1,\dots,M$, such that
\be
\ul B(\ul \alpha_j, \ul \gamma_r)=B(\alpha_j,d_r),\qquad \ul B(\ul
\gamma_r,\ul \gamma_s)=0.
\ee 

We distinguish the two different bases $\{\ul \alpha_j\}$ or
$\{\ul\alpha_j'\}$ by the subscripts $\alpha d$ or
$\alpha\gamma$. 
As a column vector in the $\alpha\gamma$ basis, the components of
${\underline k}=\sum_j k_j\ul \alpha_j+\sum_r (x_r+\nu_r) \ul \gamma_r$ read,
\be
\label{eq:ulk}
{\underline {k}}=\left(\begin{array}{c} k_1 \\ \vdots \\ k_N\\ x_1+\nu_1 \\
                       \vdots \\  x_M+\nu_M\end{array}\right)_{\!\!\!\alpha\gamma} \equiv
                   \left(\begin{array}{c} k_1 \\ \vdots \\ \\ \\
                           \vdots \\ k_{M+N} \end{array}\right) _{\!\!\!\alpha\gamma}.
                       \ee

The corresponding $(N+M)\times (N+M)$ matrix quadratic form $\underline {\bf A}$ of $\underline \Lambda$
reads
\be
\label{eq:ulA}
{\underline {\bf A}}=\left( \begin{array}{cccccc} & \vdots & & & \vdots & \\
                      \cdots & B(\alpha_i,\alpha_j)
                                                  & \cdots &  
                                                            \cdots &
                                                                     B(\alpha_i,d_r)
                      & \cdots \\ & \vdots & & &  & \\ \cdots &
                                                                B(d_s,\alpha_j)
                                                   & \cdots & & 0 &
                      \\& \vdots & & & 
                                                       &   \end{array}\right)=\left( \begin{array}{cc} {\bf A}
                  & {\bf C} \\ {\bf C}^T& 0\end{array}\right),
\ee   
with $1\leq i,j \leq N$ and $1 \leq s,r \leq M$. We introduced here the
$N\times N$ matrix ${\bf A}$ with entries $B(\alpha_i,\alpha_j)$, and
$N\times M$
matrix ${\bf C}$ with entries $B(\alpha_i,d_r)$. 

We note that the
Schur complement\footnote{Let
  $M$ be a square matrix of the form $$M=\left(\begin{array}{cc} A & B \\ C & D
                                                             \end{array}\right),$$
                                                             with
                                                             $A,B,C,D$
                                                             submatrices
                                                             of
                                                             appropriate
                                                             size. Assuming
                                                                              that
                                                                              $D$
                                                                              is
                                                                              invertible,
                                                                              the
                                                                              Schur
                                                                              complement $M/D$
                                                                              of
                                                                              $D$
                                                                              in
                                                                              $M$
                                                                              is
                                                                              the
                                                                              matrix $$M/D=A-BD^{-1}C.$$
                                                                              Moreover,
                                                                              the
                                                                              inverse
                                                                              of
                                                                              $M$
                                                                              reads
                                                                              in
                                                                              terms
                                                                              of
                                                                              the
                                                                              Schur 
                                                                              complement $M/D$, $$M^{-1}=\left(\begin{array}{cc}
                                                                                                     (M/D)^{-1}
                                                                                                     &
                                                                                                       -(M/D)^{-1}BD^{-1}
                                                                                                     \\
                                                                                                     -D^{-1}C(M/D)^{-1}
                                                                                                     & D^{-1}+D^{-1}C(M/D)^{-1}BD^{-1} \end{array}\right).$$
                                                              }
  of the block $\bf A$ in $\ul {\bf A}$ (\ref{eq:ulA}) is the matrix $\ul {\bf A}/{\bf A}=-{\bf
  C}^T{\bf A}^{-1}{\bf C}$, which will appear often below.
 We can determine the inverse of $\underline {\bf A}$ in this block form:
\be
{\underline {\bf A}}^{-1}=\left( \begin{array}{cc}  {\bf A}^{-1} -
                                   {\bf A}^{-1}{\bf C}({\bf C}^{T}{\bf A}^{-1}{\bf
                                   C})^{-1}{\bf C}^T {\bf A}^{-1}& \quad  {\bf A}^{-1}{\bf C}({\bf C}^{T}{\bf A}^{-1}{\bf
                                   C})^{-1}\\ ({\bf C}^{T}{\bf A}^{-1}{\bf
                                   C})^{-1}{\bf C}^T {\bf A}^{-1}  & \quad -({\bf C}^{T}{\bf A}^{-1}{\bf
                                   C})^{-1} \end{array}\right).
\ee 
We determine for the determinant of $\underline {\bf A}$,
\be
\det(\underline {\bf A})=\det({\bf A}) \det( -{\bf C}^T {\bf A}^{-1} {\bf
C}).
\ee 
Using that $d_s=\sum_{i,j=1}^N
B(d_s,\alpha_i) ({\bf A}^{-1})^{ij}\,\alpha_j\in \Lambda$ one
deduces that
\be
\label{DefB}
{\bf C}^T {\bf A}^{-1} {\bf
C}={\bf D},
\ee  
with ${\bf D}$ the $M\times M$ matrix with entries
$B(d_s,d_r)$. Clearly, ${\bf D}$ is positive definite since
the $d_r$ span a subspace of $\Lambda$. Using the relation
\eqref{DefB}, we deduce that $M$ independent null vectors
$\ul\gamma_r\in \ul \Lambda$, $r=1,\dots,M$, are given in terms of ${\bf A}$
and ${\bf C}$ by 
\be
\ul \gamma_r=\left( {\bf C}_{ri}^T ({\bf A}^{-1})^{ij}{\bf}\alpha_j
  ,d_r\right).
\ee
Furthermore, the
determinant of $\underline {\bf A}$ takes a simple form in terms of
the determinants of ${\bf A}$ and ${\bf D}$,
\be
\det(\underline {\bf A})=(-1)^M \det({\bf A}) \det({\bf D}).
\ee
This shows that the quadratic form $\underline {\bf A}$ is singular if the
$d_r$ are not linearly independent, in particular if $M>N$. In the following,
we will assume that the $d_r$ are linearly independent.
If $M=N$, the determinant can also be written as
\be
\det(\underline
{\bf A})=(-1)^N\det({\bf C})^2.
\ee
With a change of basis, we can bring $\underline {\bf A}$ to the block diagonal
form of the $\alpha d$ basis,
\be
\label{eq:ChofB}
{\bf G}^T {\underline {\bf A}}\, {\bf G} = \left(\begin{array}{cc}
                                                   {\bf A} & 0 \\ 0 &
                                                                      -{\bf D} \end{array}\right),
\ee
where
\be
{\bf G}=\left( \begin{array}{cc} {\bf I}_N & -{\bf A}^{-1} {\bf C}\\ 0
                                           & {\bf
                                             I}_M \end{array}\right)\in  {\rm SL}(N,M;\mathbb{Q})
                                       ,
\ee
with ${\bf I}_{\ell}$ the $\ell\times\ell$ identity matrix. We
deduce from the above that if all entries of the matrix ${\bf G} $ are integers,
thus ${\bf G}\in {\rm SL}(N,M;\mathbb{Z})$, the
lattice $\underline \Lambda$ is contained in the direct sum,
$\Lambda\oplus (-\Lambda_{d})\subset \Lambda\oplus (-\Lambda)$, with 
the lattice $\Lambda_d$ being the sublattice of $\Lambda$ generated by 
$\{d_r\}$. More generally, we have the exact sequence
\be
0\longrightarrow \Lambda \longrightarrow \ul \Lambda
\longrightarrow -\Lambda_d \longrightarrow 0.
\ee

If $M< N$, ${\bf C}^T {\bf C}$ is invertible, while ${\bf C}
{\bf C}^T$ is not. The projection to the space spanned by $\{d_r\}$, ${\bf P}:\Lambda\to \Lambda_d$,
is given by
\be
\label{eq:projection}
{\bf P}={\bf A}^{-1} {\bf C} {\bf D}^{-1} {\bf C}^T.
\ee

\subsubsection*{Elliptic Variable}
We continue with determining the elliptic variable $\underline z$ for
$\Theta_{\ul \Lambda,\ul \mu}$. We express $\ul z$ in the
$\alpha\gamma$ basis as
\be
\label{eq:zrhosigma}
\underline z=\left( \begin{array}{c} \rho \\ \sigma \end{array}
\right)_{\alpha\gamma}=\sum_{j=1}^N \rho_j \ul \alpha_j+\sum_{r=1}^M
\sigma_r \ul \gamma_r,
\ee
with $\rho$ an $N$-dimensional vector and $\sigma$ an $M$-dimensional vector. Then
\be
\ul B(\ul k, \ul z)=B(k,\rho)+B(\sum_r x_rd_r,\rho)+B(k,\sum_r \sigma_r d_r).
\ee 
Comparison with the elliptic variables in Eq. (\ref{PhiSigns}) shows that this should equal 
\be
B(k,v)+B(u,\sum_r x_r d_r).
\ee
Comparison of these two equations gives the following relations,
\be
\rho+\sum_r \sigma_r d_r=v,\qquad B(u-\rho,d_r)=0 \text{ for
  each } r=1,\dots, M.
\ee
Taking the innerproduct of the first identity with $\alpha_j$, and
some algebra using the second identity gives in terms of matrices
\be
\label{eq:sigma}
\sigma={\bf D}^{-1}{\bf C}^T (v-u)_\alpha,
\ee
where $u$ and $v$ are both vectors in the $\alpha_j$ basis, and
\be
\label{eq:rho} 
\rho=(1-{\bf P})v+{\bf P}u,
\ee
with ${\bf P}$ the projection \eqref{eq:projection}.
If we simplify to $d_j=\alpha_j$, the projection reduces to the
identity matrix, ${\bf P}={\bf 1}$, such that in this case,
\be
\label{eq:ulzMM}
\ul z=\left( \begin{array}{c} u \\ v-u \end{array} \right)_{\alpha\gamma}.
\ee

\subsubsection*{The vectors $C_r$ and $C_r'$}
We may then write Eq. (\ref{PhiSigns}) as
\begin{eqnarray}
\label{PhiIndef}
  \Phi_{\mu,\nu}(\tau, u, v,\{d_r\}) &=& \sum_{{\underline k}\in
                             \underline \Lambda + {\ul \mu}}q^{\frac{1}{2}{\underline
                                Q}(\underline k)}e^{2\pi i {\underline
                                B}(\underline k, \underline z)}\\ \nn
                             && \times 2^{-M}\prod_{r=1}^M ({\rm
                             sgn}({\underline B}(C_{r},\underline k)+\epsilon)+{\rm
                             sgn}({\underline B}(C_r', \underline
                                k+\underline a))),
\end{eqnarray}
for vectors $C_r$ and $C_r'\in \underline \Lambda^*$, and with 
\be
\label{eq:ulmu}
\ul \mu=\left(\begin{array}{c} \mu-\nu \\ \nu \end{array} \right)_{\alpha\gamma}\in
\ul \Lambda^*,
\ee 
with respect to the $\alpha\gamma$-basis (\ref{eq:ulk}). We choose
not to underline $C_r$ and $C_r'$, since these vectors do not have a
counterpart in $\Lambda$.  

We proceed by determining these vectors. Since ${\underline
  B}(C_r',{\underline k})=B(d_r,k)$ for all $\underline k$, $C_r'$
equals the basis element $\gamma_r\in \ul \Lambda$. As a vector, we have
\be
\label{eq:Crp}
{\bf C}_r'=\left(\begin{array}{c} 0 \\ \vdots \\ 0 \\1 \\ 0\\ \vdots \\
                   0 \end{array}\right)_{\alpha\gamma},
\ee
with the only non-vanishing entry the 1 on the $(N+r)'$th element. For
the norms and innerproducts of $C_r'$, we have
\be
\underline Q(C_r')=0, \qquad \underline B(C_r',C_s')=0.
\ee
This will have the consequence that the ${\rm sgn}(\cdot)$'s whose arguments
involve $C_r'$ remain unchanged in the transition from $\Phi_\mu$ to
its completion $\widehat \Phi_\mu$. See
for more details Section \ref{sec:ModCom}.
 
Furthermore since ${\underline
  B}(C_r,{\underline k})=x_r$ for all $\underline k$, $C_r$ must be
the vector in $\ul \Lambda$ dual to $C_r'$ in $\underline \Lambda$,
\be
\underline B(C_r,C_s')=\delta_{rs}.
\ee
Therefore, as a vector, ${\bf C}_r$ is given by
\be
\label{eq:Cr}
{\bf C}_r=\underline {\bf A}^{-1} {\bf C}_r'.
\ee 
The norm of $C_r$, ${\underline Q}(C_r)$, is given by the $(N+r)$'th
diagonal element of ${\underline {\bf A}}^{-1}$, or equivalently,
\be
{\underline Q}(C_r)=-({\bf C}^T {\bf A}^{-1} {\bf C})^{-1}_{rr}=-{\bf D}^{-1}_{rr}.
\ee
The innerproduct of $C_r$ and $C_s$ is given by 
\be
\label{eq:BCrCs}
{\underline B}(C_r,C_s)=-{\bf D}^{-1}_{rs}.
\ee 

It is helpful to transform ${\bf C}_r$ to the block diagonal basis
\eqref{eq:ChofB}. One then finds 
\be 
\label{eq:Cr}
{\bf G}^{-1} {\bf C}_r= -\left( \begin{array}{c} {\bf 0}_N \\
                                 {\bf D}^{-1}_{r1} \\ \vdots  \\ {\bf
                                  D}^{-1}_{rM} \end{array}
                            \right)_{\alpha d},
                             \ee
where ${\bf D}^{-1}_{ij}$ are the entries of the inverse matrix of
$\bf D$ \eqref{DefB}.
Thus the projection of $C_r\in \underline \Lambda^*$ to $\Lambda^*$ vanishes, while the projection to
$(-\Lambda_d)^*$ equals the vector $-d_r^*$, with
$d^*_r$ equal to
\be
d_r^*=\sum_{s=1}^M B(d_r,d_s)^{-1}\,d_s,
\ee 
that is to say the dual vector to $d_r$ in $\Lambda_d$, which satisfies
\be
\label{eq:dr*ds}
B(d_r^*,d_s)=\delta_{r,s}.
\ee
Thus summarizing $C_r=\ul \alpha_{N+r}^*$.

\subsubsection*{Kernel}
Finally, we need to address the fact that the argument of one of the
signs is $x_r+\varepsilon=B(C_r,\underline k)+\varepsilon$ with
$\varepsilon$ a sufficiently small positive constant, rather than $B(C_r,\underline k+\ul a)$
with ${\ul a}={\rm Im}(\ul z)/y$. Good periodicity and modular
properties require that the function can be expressed as $\ul k +\ul a$
rather than $\ul k$ and $\ul a$ separately. This may be seen for
example from the Poisson resummation technique. To achieve this, we introduce the following abbreviations
\be
\begin{split}
&{\rm s}_{r,\epsilon}=\sgn(x_r+\epsilon),\\
&{\rm s}_{r,a}=\sgn(\ul B(C_r,\underline k+\ul a))=\sgn(\nu_r+x_r+{\rm Im}(\sigma_r)/y),\\
&{\rm s}_{r,a}'=\sgn(\ul B(C_r',\underline k+\ul a))=\sgn(B(d_r,k+{\rm Im}(\rho)/y)),
\end{split}
\ee
with $\sigma_r$ determined by Eq. (\ref{eq:sigma}) and $\rho$ as in
Eq. (\ref{eq:rho}). The expression for $\rho$ demonstrates that
$B(d_r,k+{\rm Im}(\rho)/y)=B(d_r,k+{\rm Im}(u)/y)=B(d_r,k+a)$, which matches with
the terms in Eq. (\ref{PhiSigns}). 

The second line in Eq. \eqref{PhiIndef} is then written as the
kernel $K^\epsilon$,\footnote{We use that terms of the form
  $\prod_{i=1}^0 {\rm s}_{r_i,\epsilon}$ are equal to 1.}
\be
\label{eq:kernelK}
\begin{split}
K^\epsilon(\{C_r, C_r'\}, \ul k, \ul a)&=2^{-M}\prod_{r=1}^m
({\rm s}_{r,\epsilon}+{\rm s}_{r,a}')\\
&=2^{-M}\sum_{P=0}^M\sum_{\{r_1,\dots
  ,r_P,s_1,\dots ,s_{M-P}\}\in \{ 1,\dots,M\}} \prod_{i=1}^P
{\rm s}_{r_i,\epsilon}\, \prod_{j=1}^{M-P} {\rm s}_{r_j,a}'.
\end{split}
\ee
where $\{r_1,\dots, r_P\}$ and $\{s_1,\dots, s_{M-P}\}$ are 
an unordered $P$ and $(M-P)$-tuple respectively. Thus in terms of
$K^\epsilon$, $\Phi_{\mu,\nu}$ reads
\be
\label{eq:PhiKernel}
\Phi_{\mu,\nu}(\tau,u,v,\{d_r\})=\sum_{\ul k\in \ul \Lambda+\ul \mu} K^\epsilon(\{C_r, C_r'\}, \ul k, \ul
a)\, q^{\frac{1}{2}{\ul Q}(\ul k)+\ul B(\ul k, \ul z)}.
\ee

To prepare for the formulation of the modular completion in the next
section, we aim to replace the
${\rm s}_{r,\epsilon}$ with ${\rm s}_{r,a}$. We therefore express the
kernel $K^\epsilon$ as the kernel $K$ defined in
Eq. (\ref{eq:KIndefTheta}) plus a term depending on $\epsilon$,
\be
\label{KCCr2}
\begin{split}
&K^\epsilon(\{C_r,C_r'\}, \underline k, \ul a)=2^{-M}\sum_{P=0}^M\sum_{\{r_1,\dots
  r_P,s_1,\dots ,s_{M-P}\}\in \{ 1,\dots,M\}} \prod_{i=1}^P
 {\rm s}_{r_i,a}\,  \prod_{j=1}^{M-P}  {\rm s}_{s_j,a}' \\
&\quad  + 2^{-M}\sum_{P=1}^{M}  \sum_{\{r_1,\dots
  r_P,s_1,\dots ,s_{M-P}\}\in \{ 1,\dots,M\}}\left(\prod_{i=1}^P
{\rm s}_{r_i,\epsilon}-\prod_{i=1}^P
{\rm s}_{r_i,a}\right) \prod_{j=1}^{M-P}{\rm s}_{s_j,a}'\\
&=K(\{C_r,C_r'\}, \underline k+\ul a)\\
&\quad + 2^{-M}\sum_{P=1}^{M}  \sum_{\{r_1,\dots
  r_P,s_1,\dots ,s_{M-P}\}\in \{ 1,\dots,M\}}\left(\prod_{i=1}^P
{\rm s}_{r_i,\epsilon}-\prod_{i=1}^P
{\rm s}_{r_i,a}\right) \prod_{j=1}^{M-P}{\rm s}_{s_j,a}'.
\end{split}
\ee
We define now the function $S_{\mu,\nu}$ as the series whose kernel is
given by minus the terms on the second line of Eq. \eqref{KCCr2},
\be
\begin{split}
S_{\mu,\nu}(\tau,u,v,\{d_r\})&=2^{-M}\sum_{k\in \Lambda+\mu-\nu\atop x_r\in \mathbb{Z}} \sum_{P=1}^{M}  \sum_{\{r_1,\dots
  r_P,s_1,\dots ,s_{M-P}\}\atop\in \{ 1,\dots,M\}} \left( \prod_{i=1}^P
  {\rm s}_{r_i,a}-\prod_{i=1}^P
  {\rm s}_{r_i,\epsilon}\, \right)
\prod_{j=1}^{M-P}{\rm s}_{s_j,a}'\\
&\times \,q^{Q(k)/2+B(\nu+x_rd_r,k)}e^{2\pi i B(u,\nu+x_r d_r)+2\pi iB(v,k)}.
\end{split} 
\ee 
We define moreover $\Phi_{\mu,\nu}^+$ as the function whose kernel is
$K(\{C_r,C_r'\}, \underline k+\ul a)$. The three functions are thus related as
\be
\label{eq:DefPhi+}
\Phi^+_{\mu,\nu}(\tau,u,v,\{d_r\})=\Phi_{\mu,\nu}(\tau,u,v,\{d_r\})+S_{\mu,\nu}(\tau,u,v,\{d_r\}).
\ee
One of the special properties of $\Phi^+_{\mu,\nu}$ is that it is periodic in both $\mu$ and $\nu$,
\be
\label{eq:Phi+period}
\Phi^+_{\mu+m,\nu+n}(\tau,u,v,\{d_r\})=\Phi^+_{\mu,\nu}(\tau,u,v,\{d_r\}),\qquad
m\in \Lambda, n\in \Lambda_d,
\ee
while $\Phi_{\mu,\nu}$ and $S_{\mu,\nu}$ do not satisfy the
periodicity in $\nu$ separately. See
Eq. (\ref{eq:nuqperiod}). Note that depending on $a$ and $\nu$, $S_{\mu,\nu}$
may vanish. We can carry out the geometric sums of
$\Phi^+_{\mu,\nu}(\tau,u,v,\{d_r\})$,
\be
\label{eq:Phi+Appell}
\begin{split}
\Phi^+_{\mu,\nu}(\tau,u,v,\{d_r\})&=e^{2\pi i B(\nu-\lfloor \nu_r + {\rm
    Im}(\sigma_r)/y\rfloor d_r,u)}\\
&\times \sum_{k\in \Lambda+\mu-\nu}
\frac{q^{\frac{1}{2}Q(k)+B(\nu-\lfloor \nu_r + {\rm
      Im}(\sigma_r)/y\rfloor d_r,k) }e^{2\pi i B(v,k)}}{\prod_{r=1}^M
  \left(1-e^{2\pi i B(d_r,u)}q^{B(d_r,k)}\right)},
\end{split}
\ee 
where $\nu_r$ are the coefficients of $\nu=\sum_r \nu_r d_r$, and
similarly for $\sigma=\sum_r \sigma_rd_r$. We have thus the relation 
\be
\label{eq:PhiPlusPhi}
\Phi^+_{\mu,\nu}(\tau,u,v,\{d_r\})=\Phi_{\mu,\tilde \nu}(\tau,u,v,\{d_r\}),
\ee
with
\be
\tilde \nu=\nu-\sum_{r=1}^M \lfloor \nu_r + {\rm
    Im}(\sigma_r)/y\rfloor d_r.
\ee

It is the function $\Phi^+_{\mu,\nu}$, which naturally takes the
form of an indefinite theta series $\Theta_{\ul \Lambda,\ul \mu}$ (\ref{eq:DefIndefTheta})
\be
\label{eq:Phi+Theta}
\Phi^+_{\mu,\nu}(\tau,u,v,\{d_r\})=\Theta_{\ul \Lambda,\ul \mu}(\tau,\ul
z,\{C_r,C_r'\}),
\ee
with $\ul \Lambda$ as described above, $\ul \mu$ by
Eq. \eqref{eq:ulmu}, $\ul z=(\rho,\sigma)$ defined by Eqs
(\ref{eq:sigma}) and (\ref{eq:rho}), and ${\bf C}_r$, ${\bf C}_r'$
defined by Eqs (\ref{eq:Crp}) and (\ref{eq:Cr}).
\vspace{.3cm}\\
{\bf Example}\\
We conclude this subsection with an elementary example illustrating the
characteristic properties of $\Phi_{\mu,\nu}$, $\Phi^+_{\mu,\nu}$ and $S_{\mu,\nu}$. For
$\mu,\nu\notin \mathbb{Z}$, we define
\be
\phi^+_{\mu,\nu}(\tau)=\tfrac{1}{2}\sum_{k \in \mathbb{Z}+\mu\atop \ell \in
  \mathbb{Z}+\nu} \left( \sgn(\ell)+\sgn(k)\right)\,q^{k^2/2+k\ell}.
\ee
The sum is clearly periodic in both $\mu,\nu $ since
this just shifts the sum over $k$ and $\ell$. The sum over $\ell$ can be
done as a geometric sum,
\be
\phi^+_{\mu,\nu}(\tau)=\sum_{k\in \mathbb{Z}+\mu} \frac{q^{k^2/2+k \{\nu\}}}{1-q^{k}},
\ee
which is also clearly periodic in $\mu,\nu$. The decomposition in the
$\phi_{\mu,\nu}$ and $s_{\mu,\nu}$ corresponds to
\be
\begin{split} 
\phi_{\mu,\nu}(\tau)=\sum_{k\in \mathbb{Z}+\mu} \frac{q^{k^2/2+k\nu}}{1-q^{k}},\qquad
s_{\mu,\nu}(\tau)=\sum_{k\in \mathbb{Z}+\mu} \frac{q^{k^2/2}}{1-q^{k}}\,(q^{k \{\nu\}}-q^{k \nu}).
\end{split}
\ee
The function $s_{\mu,\nu}$ can be written as
\be
s_{\mu,\nu}(\tau)=\left(\sum_{k\in \mathbb{Z}+\mu+\nu} q^{k^2/2}
\right)\times \left\{\begin{array}{cc}  \sum_{m=0}^{\lfloor \nu \rfloor -1}
                  q^{-(m+\{\nu\})^2/2},&\qquad \nu>0,\\ \sum_{m=0}^{-\lfloor \nu \rfloor -1}
                  q^{-(m+\{-\nu\})^2/2}, &\qquad \nu<0.\end{array}\right.
\ee
Thus $s_{\mu,\nu}$ is a theta series times a finite number of terms,
ie $|\nu-\{\nu\}|=|\lfloor \nu \rfloor|$ terms, which vanishes for $\nu=\{\nu\}$.

\subsection{Specialization to the Root Lattice $A_N$ and its Weyl Reflections}
\label{sec:Anlattice}
In this section we choose the lattice $\Lambda$ as the root lattice $A_N$
of
the Lie group $SU(N+1)$. The dimension of $\Lambda$ equals the rank
$N$ of $SU(N+1)$, and we choose for the generators $\alpha_j$ of
$\Lambda$ the simple roots $SU(N+1)$.
Then $d_r=\sum_{j=1}^N
d_{r,j} \alpha_j$ where, $d_{r,j}\in\mathbb{Z}$. If we choose
$d_j=\alpha_j$ for $j=1,\dots, N$, we shorten the notation of $\Phi_{\mu}$ as follows:
$$\Phi_{\mu,\nu}(\tau,u,v,\{ \alpha_j \})=\Phi_{\mu,\nu}(\tau,u,v).$$

The matrices ${\bf D}$ and ${\bf C}$ of Section
\ref{sec:AppellIndef} are then both equal to $SU(N+1)$ Cartan the matrix ${\bf A}$, and ${\bf P}={\bf
  1}$. The matrix quadratic form then reads
\be
\ul{\bf{A}}=\begin{pmatrix}
{\bf A} & {\bf A}\\
{\bf A} & 0
\end{pmatrix}, \quad \ul{\bf A}^{-1}=\begin{pmatrix}
0 & {\bf A}^{-1}\\
{\bf A}^{-1} & -{\bf A}^{-1}
\end{pmatrix}.
\ee
Using the inverse of the Cartan matrix ${\bf A}$ of $SU(N+1)$ \footnote{The form can be checked easily as right (left) inverse by writing the $i$-th row (column) of ${\bf A}_{ij}= 2\delta_{i,j}-\delta_{i,j-1}-\delta_{i,j+1}$ when $1<i<N$ for $i\in\{1,\dots N\}$ writing ${\bf A}_{1j}= 2\delta_{1,j}-\delta_{1,j-1}$, ${\bf A}_{Nj}=2\delta_{N,j}-\delta_{N,j+1}$ and taking the inner product with ${\bf A}^{-1}_{jl}={\rm min}(j,l)-\frac{jl}{N+1}$ and observing $l\ge k$ and $l<k$ cases separately.} we have for the inner products
\be
\begin{split}
& {\underline Q}(C_j)=\frac{j^2}{N+1}-j<0,\qquad {\underline
  B}(C_i,C_j)=\frac{ij}{N+1}-{\rm min}(i,j),\\
& {\underline
  B}(C_i,C_j')=\delta_{i,j},\qquad \underline B(C_i',C_j')=0.
\end{split}
\ee

For later sections, it will be useful to consider a set $d_j$, which are
related to the simple roots $\{\alpha_j\}$ by the Weyl reflection
$S_{\alpha_m}$ with respect to $\alpha_m$ for some $m$, thus $d_j=S_{\alpha_m}(\alpha_j)$. 
We will look at how the Weyl reflection affects $\Phi_{\mu,\nu}$ ie, when $d_r=S_{\alpha_m}(\alpha_j)$.
See Appendix \ref{AppWeylSym} for the precise definition. From the action of the Weyl transformation $S_{\alpha_m}$ on the roots
$\alpha_j$ (\ref{weylnj}), we deduce that the components $k'_j$ of
$k'=S_{\alpha_m}(k)$ with $k\in \Lambda$ are 
\be
k'_m=-k_m+k_{m+1}+k_{m-1}, \quad k'_j=k_j \; {\rm if}\; j\ne m,
\ee
or with the matrix entries
\be
\label{eq:Salpha}
\begin{split}
&(S_{\alpha_m})_{m,m}= -1,\\
&(S_{\alpha_m})_{m,m\pm 1}=1,\\
&(S_{\alpha_m})_{m\pm 1,m}=0,\\
&(S_{\alpha_m})_{ij}=\delta_{i,j},\quad {\rm if}\quad  i\ne m,  m\pm
1, \vee \,\, j\ne m,  m\pm 1.\\
\end{split}
\ee 
One easily verifies that $S_{\alpha_m}=S_{\alpha_m}^{-1}$. Moreover,
$S_{\alpha_m}$ leaves the lattice $\Lambda$ and its quadratic 
form invariant as in Eq. (\ref{eq:TrafoInv}), such that (\ref{eq:PhiGInv}) holds
for $G=S_{\alpha_m}$.

The $A_N$ lattice has $N+1$ conjugacy classes $\mu\in \Lambda^*/\Lambda$. 
The Weyl group leaves the conjugacy class specified by
$\mu=(\mu_1,\dots, \mu_N)$ invariant,
\be
\label{mu'mu}
\mu_m'=\mu_{m+1}+\mu_{m-1}-\mu_{m}=\mu_m \mod \mathbb{Z},
\ee
such that $\Lambda+\mu$ remains invariant. In this way, we can
determine the $N+1$ conjugacy classes
algorithmically. The second equality of Eq. (\ref{mu'mu}) gives the following relations,
  \begin{eqnarray}
2\mu_1 &=& \mu_2,\\ \nn
2\mu_2 &=& \mu_1+\mu_3 \mod \mathbb{Z},\\ \nn
2\mu_3 &=& \mu_2+\mu_4 \mod \mathbb{Z},\\ \nn
&\vdots &\\ \nn
2\mu_N &=& \mu_{N-1}\mod \mathbb{Z},
\end{eqnarray}
which can be solved by
\begin{eqnarray}\label{mu_res}
\mu_j=j\mu_1\mod \mathbb{Z}, \quad (N+1)\mu_1=0.
\end{eqnarray}
The above implies that
if we restrict $0\le \mu_1 <1$, then 
\be
\mu_1+\mu_{N}\in\{0,1\},\quad 0 \le \mu_j<1.
\ee
This gives for the conjugacy classes for $N=2,3$:
\begin{eqnarray}
A_2 &:& (\mu_1,\mu_2)=(0,0),\;\; (\tfrac{1}{3},\tfrac{2}{3}),\;\; (\tfrac{2}{3},\tfrac{1}{3}),\\ \nn
A_3 &:& (\mu_1,\mu_2,\mu_3)=(0,0,0),\;\;(\tfrac{1}{4},\tfrac{1}{2},\tfrac{3}{4}),\;\;(\tfrac{1}{2},0,\tfrac{1}{2}),\;\; (\tfrac{3}{4},\tfrac{1}{2},\tfrac{1}{4}).
\end{eqnarray}

In many cases it is useful to mod out the set of conjugacy classes by
the $\mathbb{Z}_2$ transformation  $\mu\to -\mu$ (\ref{Phi-mu}), this
gives ${\frac{N+3}{2}}$ inequivalent classes for $N$ odd and
${\frac{N+2}{2}}$ for $N$ even. The set of conjugacy classes modulo
this action can obtained by restricting $0\le \mu_1 \le \frac{1}{2}$
with 
\be
\mu_1+\mu_{N}=0\mod \mathbb{Z}.
\ee
This gives for $N=2,3$:
\begin{eqnarray}
A_2 &:& (\mu_1,\mu_2)=(0,0),\;\; (\tfrac{1}{3},-\tfrac{1}{3}),\\ \nn
A_3 &:& (\mu_1,\mu_2,\mu_3)=(0,0,0),\;\;(\tfrac{1}{4},\tfrac{1}{2},-\tfrac{1}{4}),\;\;(\tfrac{1}{2},0,-\tfrac{1}{2}).
\end{eqnarray}

While $S_{\alpha_m}$ leaves invariant ${\bf A}$, it does affect the
extended lattice $\ul{\mathbf{A}}$.
The explicit form of $\ul{\mathbf{A}}_m$ is given in terms of a
deviation from $\ul{\mathbf{A}}$ as:
\be\ul{\mathbf{A}}_m=\begin{pmatrix}
\mathbf{A} & S_{\alpha_m}^T \mathbf{A}\\
S_{\alpha_m}^T \mathbf{A} & 0
\end{pmatrix}.
\ee
with inverse,
 \begin{eqnarray}
 \ul{\mathbf{A}}_m^{-1} &=& \left(\begin{array}{cc}
 \mathbf{0} & {\bf A}^{-1}-{\bf 1}_{mm}\\
 {\bf A}^{-1}-{\bf 1}_{mm} & -\bf A^{-1}
 \end{array}\right),
 \end{eqnarray}
where the matrix ${\bf 1}_{mm}$ is the $N\times N$ matrix with only
the $(mm)$-th entry equal to 1 and the rest zero.

\section{Modular Completion of Appell Functions}
\label{sec:ModCom}
We determine in this section the modular completion of the Appell
functions using the relation to indefinite theta series derived in the
previous section. The modular completion of indefinite theta series
are well-established \cite{ZwegersThesis, Alexandrov:2016enp, Nazaroglu:2016lmr} using the results of Vign\'eras \cite{Vigneras:1977}.
We will first consider the modular completion $\widehat \Phi_\mu$ of $\Phi_\mu$
(\ref{phigen1}), since the completion of $\Phi_{\mu,\nu}$ can be
derived from $\widehat \Phi_\mu$ using Eq. (\ref{eq:PhimunuPhimu0}). Our main results are the
structural formulas (\ref{eq:StructForm}) for $\widehat \Phi_\mu$ and
(\ref{eq:StructForm2}) for $\widehat \Phi_{\mu,\nu}$,
which demonstrates that these involve functions
$\Phi_{\mu',\nu'}(\tau,u,v,\{d_r\}')$, with $\{d_r\}'$ subsets of
$\{d_r\}$. 

The resulting functions $\widehat \Phi_\mu$ and $\widehat \Phi_{\mu,\nu}$ transform as
\be
\begin{split}
&\widehat \Phi_{\mu,\nu}\left(\frac{a\tau+b}{c\tau+d}, \frac{a\bar
    \tau+b}{c\bar \tau+d},
  \frac{u}{c\tau+d}, \frac{\bar u}{c\bar \tau+d}, \frac{v}{c\tau+d},
  \frac{\bar v}{c\bar \tau+d};\{d_r\} \right)\\
&\quad  =
(c\tau+d)^{(M+N)/2}  \exp\left(-\pi i \frac{c}{c\tau+d} \left(Q(u^2)
    -2B(u,v)\right)\right)\\
&\qquad \times \widehat \Phi_{\mu,\nu}(\tau,\bar \tau,u,\bar u,v,\bar v;\{d_r\}),
\end{split}
\ee 
for $\left(\begin{array}{cc} a & b \\ c & d \end{array}\right)\in
\Gamma(4n)$ with $n=|{\rm det}(\Lambda)\,{\rm
  \det}(\Lambda_d)|$. Depending on details of the lattice, $n$ maybe
smaller than this value.

\subsection{Modular Completion of Indefinite Theta Series}
We consider first the indefinite theta series $\Theta_{\Gamma,  \mu}$ as introduced in
Eq. (\ref{eq:DefIndefTheta}). This function does not transform as a
modular form, essentially because the sum does not have support on the
full lattice $\Gamma$ (or a sublattice). We recall now the
modular completion $\widehat \Theta_{\Gamma, \mu}$ of $\Theta_{\Gamma,
  \mu}$ as in \cite{Alexandrov:2016enp, Nazaroglu:2016lmr}. 

For an indefinite lattice $\Gamma$, we recall the definition of the generalized
error function $\CE_P$ and complementary generalized error function
$\CM_P$ \cite[Eq. (6.1) and (6.3)]{Alexandrov:2016enp}\footnote{Note
  both $\CE_P$ and $\CM_P$ differ by the factor $(-1)^P$ due to the
  opposite sign for the convention of the quadratic form.}
\be 
\label{eq:DefEPMP}
\begin{split}
&\CE_P(\{C_j\}, \Gamma; x)= (-1)^P\int_{\langle \{ C_j\}\rangle}
\prod_{j=1}^{P} {\rm sgn}( B(C_j,y))\,e^{\pi  Q(y- x^{||})}\, d^{P}y,\\
&\CM_{P}(\{C_j\}, \Gamma ; x) =
\sqrt{|\Delta(\{C_j^{\star}\})|}\left(\frac{1}{\pi i}\right)^{P}\\
&\qquad \times \int_{\langle
 \{C_j\}\rangle-ix^{||}}\prod_{j=1}^{P} \frac{1}{ B(C_j^\star,z)}
e^{\pi  Q(z)+2\pi i  B(z,x)} d^{P}z,
\end{split}
\ee 
where $\{C_j\}=C_1,\dots, C_{P}$ are independent time-like vectors,
$Q(C_j)<0$, spanning $P$ directions, $x\in  \Gamma\otimes
\mathbb{R}$ and $x^{||}$ is the orthogonal projection to the plane
spanned by $\{C_j\}$. $\{C_j^\star\}$ is the dual basis to
$\{C_j\}$ in the plane spanned by $\{C_j\}$.\footnote{Note the subtle
  distinction we make between $*$ and $\star$. $\{C_j^*\}$ is a subset of
  vectors which are dual to $C_r$, for all $r$, in the full lattice $\Lambda$, while $\{C_j^\star\}$ is the set of dual
vectors to the set $\{C_j\}$ in the plane spanned by the set $\{C_j\}$.} Moreover, $\Delta(\{C_j^{\star}\})$ is the determinant of the Gram
matrix $B(C_i^\star, C_j^\star)$. The domain of the $P$-dimensional integral is the $P$-plane
spanned by $\{C_j\}$, and is normalized such that
\be 
\int_{\langle \{C_j\}\rangle} e^{\pi Q(y)}d^{P}y=1.
\ee  

Nazaroglu \cite[Prop. 3.15]{Nazaroglu:2016lmr} expresses the $\CE_P$ in terms of
the functions $\CM_P$. In our notation this reads, 
\be
\label{eq:EPMP} 
\begin{split}
\CE_P(\{C_r\},x; \Gamma)&=\sum_{L=0}^P \sum_{\{ v_1,\dots, v_L, w_{1},\dots ,
  w_{P-L}  \}\in \{1,\dots,P\}}     \CM_L( \{C_{v_i}\}, \Gamma; x) \\
& \times \prod_{j=1}^{P-L} \sgn( B(C^{\perp V_L}_{w_j},x)),
\end{split}
\ee 
where $V_L$ is the hyperplane spanned by the $\{C_{v_i}\}$, and
$C^{\perp V_L}_w$ is the component of $C_w$ orthogonal to $V_L$.

The completion $\widehat \Theta_{ \Gamma,  \mu}(\tau, z,
\{C_r,C_r'\})$ is obtained from $\Theta_{\Gamma,  \mu}(\tau, z,
\{C_r,C_r'\})$ by replacing all products of signs by generalized error
functions. That is to say, the completion is defined as
\be
\widehat \Theta_{ \Gamma,  \mu}(\tau, \bar \tau, z, \bar z,
\{C_r,C_r'\})=\sum_{k\in \Gamma+\mu}
\widehat K(\{C_r,C_r'\}, k+  a)\,q^{ Q(k)/2}e^{2\pi i  B( k, z)},
\ee 
with
\be
\widehat K(\{C_r,C_r'\}, k)=2^{-M}\sum_{P=1}^{M}  \sum_{\{v_1,\dots
  v_P,s_1,\dots ,s_{M-P}\}\in \{ 1,\dots,M\}} \CE_M(\{C_v,C'_s\},\Gamma;\sqrt{2y}\, k).
\ee

\subsection{Kernels}

To determine the modular completion of $\Phi_\mu$, we recall
Eq. (\ref{eq:Phi+Theta}). We thus define the modular completion
$\widehat\Phi_\mu$ of $\Phi_\mu$ as the modular completion $\widehat
\Theta_{\ul \Lambda,\ul \mu}$. We thus replace the  
the kernel $K^\epsilon$ by $\widehat K$. Since the vectors $C_r'$ have
a vanishing norm, the corresponding $\sgn$'s are not modified in the
completion. Therefore,  $\widehat K$ is obtained from $K^\epsilon$ by replacing
$\prod_{i=1}^P{\rm s}_{r_i,a}$ by $\CE_P(\{C_{r_i}\},x;\underline \Lambda)$ and subtracting
the second line in Eq. \eqref{KCCr2}. With $y={\rm Im}(\tau)$, the kernel
$\widehat K$ for the completed function thus reads
\be 
\begin{split}
  &\widehat K(\{C_r,C_r'\}, \underline k, \underline
  a)=K^\epsilon(\{C_r,C_r'\}, \underline k, \underline a) +  2^{-M}\sum_{P=1}^{M}  \sum_{\{r_1,\dots
  r_P,s_1,\dots ,s_{M-P}\}\in \{ 1,\dots,M\}}\\
&\qquad  \left( \CE_P(\{C_{r_i}\},
  \sqrt{2y}(\underline k+\underline a))-\prod_{i=1}^P
  {\rm s}_{r_i,\epsilon}\, \right) \prod_{j=1}^{M-P}{\rm s}_{s_j,a}' .
\end{split}
\ee
This equals
\be
\label{eq:CompKernel}
\begin{split}
  &\widehat K(\{C_r,C_r'\}, \underline k, \underline
  a)=K^\epsilon(\{C_r,C_r'\}, \underline k, \underline a) \\
&+ 2^{-M}\sum_{P=1}^{M}  \sum_{\{r_1,\dots
  r_P,s_1,\dots ,s_{M-P}\}\atop\in \{ 1,\dots,M\}} \left( \prod_{i=1}^P
  {\rm s}_{r_i,a}-\prod_{i=1}^P 
  {\rm s}_{r_i,\epsilon}\, \right) \prod_{j=1}^{M-P}{\rm s}_{s_j,a}' \\
  &+2^{-M}\sum_{P=1}^{M}  \sum_{\{r_1,\dots
  r_P,s_1,\dots ,s_{M-P}\}\atop \in \{ 1,\dots,M\}} \left( \CE_P(\{C_{r_i}\}, \underline \Lambda;
  \sqrt{2y}(\underline k+\underline a))-\prod_{i=1}^P
  {\rm s}_{r_i,a}\, \right) \prod_{j=1}^{M-P}{\rm s}_{s_j,a}' 
\end{split}
\ee
As discussed above, the lattice for the Appell functions splits, $\underline \Lambda=\Lambda \oplus (-\Lambda_d)$ and
the $C_r$ are orthogonal to the positive definite lattice
$\Lambda$. As a result, the
function $\CE_P$ \eqref{eq:DefEPMP} simplies as we discuss later.

\subsection{Modular Completion}
The first term on the
rhs gives $\Phi_\mu(\tau,u,v,\{d_r\})$. The second line
gives further holomorphic terms, which can vanish in many cases, for
example for certain non-vanishing $\nu$. The third line in
\eqref{eq:CompKernel} is non-holomorphic and vanishes in the limit
$y\to\infty$ (assuming that $\ul k+\ul a$ is non-vanishing). This term is our main interest.

The modular completion $\widehat \Phi_\mu$ of $\Phi_\mu$ is obtained
by replacing the kernel $K^\epsilon$ in Eq. (\ref{eq:PhiKernel}) by
$\widehat K$,
\be
\label{eq:PhiHat}
\begin{split}
&\widehat \Phi_\mu(\tau,\bar \tau, u,\bar u, v, \bar v, \{d_r\})= \sum_{x_r\in
  \mathbb{Z}} \sum_{k\in \Lambda+\mu} \widehat
K(\{C_r,C_r'\},\underline k,\underline a )\\
&\quad \times q^{Q(k)/2+B(x_rd_r,k)}e^{2\pi i B(v,k)+2\pi i B(x_rd_r
  u)}.
\end{split}
\ee
In the following, we will
\begin{itemize}
\item split holomorphic and non-holomorphic terms,
\item write the non-holomorphic terms in terms of data associated to the
  lattice $\Lambda$ rather than $\ul \Lambda$,
\item write the non-holomorphic terms in terms of $\Phi_{\mu,\nu}^+$ with a
  smaller depth $M'<M$.
\end{itemize}

In terms of $\Phi^+_\mu$ (\ref{eq:DefPhi+}), we split the holomorphic
and non-holomorphic parts of $\widehat \Phi_\mu$,
\be
\begin{split}
&\widehat \Phi_\mu(\tau,\bar \tau, u,\bar u, v, \bar v,\{d_r\})=
\Phi^+_\mu(\tau,u,v,\{d_r\}) + \CR_\mu(\tau,\bar \tau, u,\bar u, v, \bar v,\{d_r\}),
\end{split}
\ee
where the non-holomorphic part $\CR$ is defined as
\be
\begin{split}
&\CR_\mu(\tau,\bar \tau, u,\bar u, v, \bar v,\{d_r\})= 2^{-M}\sum_{P=1}^{M}  \sum_{\{r_1,\dots
  r_P,s_1,\dots ,s_{M-P}\}\in \{ 1,\dots,M\}}\\
&\qquad   \sum_{x_{r_i}\in \mathbb{Z}} \left( \CE_P(\{C_{r_i}\}, \ul \Lambda;
  \sqrt{2y}(\underline k+\underline a))-\prod_{i=1}^P
  {\rm s}_{r_i,a}\, \right) \\
&\qquad   \sum_{k\in \Lambda+\mu} \sum_{x_{s_\ell}\in \mathbb{Z}} \left(\prod_{\ell=1}^{M-P}{\rm s}_{s_\ell,a}'\right)
q^{Q(k)/2+B(x_rd_r,k)}e^{2\pi i B(v,k)+2\pi i B(x_rd_r,u)}.
\end{split}
\ee

We continue by expressing the $\CE_P$ in terms of $\CM_P$ using
Eq. (\ref{eq:EPMP}). We thus arrive at
\be
\begin{split}
&\CR_\mu(\tau,\bar \tau, u,\bar u, v, \bar v,\{d_r\})=2^{-M} \sum_{P=1}^{M} \sum_{L=1}^{P} \sum_{\{v_1,\dots,v_L,w_{1},\dots,
  w_{P-L},s_1,\dots ,s_{M-P}\}\in \{ 1,\dots,M\}}\\
&\qquad   \sum_{x_{r_i}\in \mathbb{Z}}  \CM_L(\{C_{v_i}\}, \ul \Lambda;
  \sqrt{2y}(\underline k+\underline a))\prod_{j=1}^{P-L}
  \sgn(B(C_{w_j}^{\perp V_L},\underline k+\underline a)) \\
&\qquad   \sum_{k\in \Lambda+\mu} \sum_{x_{s_\ell}\in \mathbb{Z}} \left(\prod_{\ell=1}^{M-P}{\rm s}_{s_\ell,a}'\right)
q^{Q(k)/2+B(x_rd_r,k)}e^{2\pi i B(v,k)+2\pi i B(x_rd_r,u)}.
\end{split}
\ee 
For the sum over $w_i$ and $s_j$, we substitute the kernel
Eq. (\ref{eq:KIndefTheta}),
\be
\label{eq:PhiHatML}
 \begin{split}
 &\CR_\mu(\tau,\bar \tau, u,\bar u, v, \bar v,\{d_r\})=
  \sum_{L=1}^{M} \sum_{\{v_1,\dots,v_L,s_{1},\dots,
   s_{M-L}\}\in \{ 1,\dots,M\}}\\
&\qquad   \sum_{k\in \Lambda+\mu}  \sum_{x_{r}\in \mathbb{Z}} 2^{-L} \CM_L(\{C_{v_i}\}, 
  \sqrt{2y}(\underline k+\underline
  a);\ul \Lambda)\,K(\{C_s^{\perp V_L},C_{s}'\}_{M-L},\underline k+\underline
  a)  \\ 
&\qquad  
\times q^{Q(k)/2+B(x_rd_r,k)}e^{2\pi i B(v,k)+2\pi i B(x_rd_r,u)}.
 \end{split}
\ee

We will now demonstrate that under suitable identification, the sum over the $k\in \Lambda+\mu$ and
$x_{s_j}=k_{N+s_j}\in \mathbb{Z}$ will combine to $\Phi^+_{\mu,\nu}$
with a subset of the vectors $\{d_r\}$. To this end, we introduce the $L$-dimensional sublattice $\underline \Lambda(\{C_{v_i}\})
\subset \underline \Lambda$ spanned by the $C_{v_i}$, $i=1,\dots,L$. We then
decompose the other vectors $C_{s_j}$, $j=1,\dots,M-L$, as
\be
C_{s}=C_{s}^{||V_L}+C_{s}^{\perp V_L},
\ee 
where $C_{s}^{||}$ and $C_{s}^{\perp}$  are the components
of $C_{s}$ parallel and orthogonal to
$\underline \Lambda(\{C_{v_i}\})$. We expand the $C_{s}^{||V_L}$ in
terms of coefficients $c_{s_i,v_j}$ as 
\be
\label{eq:CsiPar}
C_{s_{i}}^{||V_L}= \sum_{j=1}^{L} c_{s_i,v_j} C_{v_j},
\ee
To specify the coefficients $c_{s_i,v_j}$, let $\underline {\tilde B}$ be the $L\times L$ submatrix of $B$ defined by the entries
\be
\underline {\tilde B}_{v_j,v_k}=\underline B(C_{v_j},C_{v_k}),\qquad j,k=1,\dots,L.
\ee 
The coefficients $c_{s_i,v_j}$ are then given by
\be
\label{eq:csv}
c_{s_i,v_j}=\sum_{k=1}^L \underline {\tilde B}_{v_j,v_k}^{-1}\, \underline
B(C_{v_k},C_{s_i}).
\ee
This can be completely expressed in terms of the quadratic form of the
lattice $\Lambda_d$ using Eq. \eqref{eq:BCrCs}. With $v_j$, $s_j$ as in the sum
(\ref{eq:PhiHatML}), we introduce the following sets $R$, $V$ and $R/V$,
\be
\begin{split}
& R=\{1,\dots, M\},\\
& V\subseteq R,\qquad V=\{v_j\in R\,|\, j=1,\dots, L\},\\
&R/V \subseteq R,\qquad R/V=\{s_j\in R \,|\,j=1,\dots, M-L\}.
\end{split}
\ee

We let ${\bf D}_V$ be the submatrix of ${\bf D}$ defined below
Eq. (\ref{DefB}), with indices in the
set $V\times V$. This shows that the matrix with entries $\tilde
{\ul B}^{-1}_{v_j,v_k}$ is the negative of the Schur complement $\tilde {\bf D}_V$ of ${\bf
  D}_{R/V}$ in the full matrix ${\bf D}$. 
We thus have
\be
\tilde {\bf D}_V={\bf D}/{\bf D}_{R/V}.
\ee 

From the definition of $\CM_L$ in Eq. (\ref{eq:DefEPMP}), we deduce that
this function in Eq. (\ref{eq:PhiHatML}) only depends on the component 
of $\ul k$ parallel to $\ul \Lambda(\{C_{v_i}\})$.
Concretely, the component of $\ul k$ parallel to the lattice $\ul \Lambda(\{
C_{v_i}\})$, $\ul k^{||V}$, is given by
\be
\begin{split}
\ul k^{||V}&=\sum_{j,k=1}^L \tilde {\ul B}^{-1}_{v_j,v_k}\, \ul
B(C_{v_k},\ul k)\,C_{v_j}\\
& =\sum_{j,k=1}^L \tilde {\ul B}^{-1}_{v_j,v_k}\, x_{v_k} C_{v_j},
\end{split}
\ee
such that 
\be
\begin{split}
\ul k^{||V}=(0, (\sum_{j,k} \tilde {\bf D}_V)_{k,j} \,x_{v_k}\,d_{v_j}^*)\in \Lambda \oplus (-\Lambda_d)^*.
\end{split}
\ee
We thus confirm that $\CM_L(\{C_{v_i}\}, \ul \Lambda;
  \sqrt{2y}(\underline k+\underline
  a))$ only depends on $x_{v_\ell}$, not on the $x_{s_i}$. As a result,  $\CE_P$
  and $\CM_P$ simplify and can be expressed in terms of the positive
  definite lattice $\Lambda$. If we parametrize the integrand by
\be
z=\sum_{i} z_i C_{r_i},
\ee 
it follows by Eq. (\ref{eq:Cr}) that
\be
\CE_P( \{C_r\},\underline \Lambda;\ul x)=E_P(\{d^*_r\},\Lambda ; x),
\ee 
with
\be
E_P(\{c_r\},\Lambda; x)=\int_{\left< \{c_r\}\right>}
\prod_{j=1}^P \sgn(B(c_j,x))\,e^{-\pi Q(y-x^{||})}\,d^Py,
\ee
and $x^{||}$ is as before the orthogonal projection of $x$ to the plane
spanned by $\{c_r\}$. 

Similarly, $\CM_P$ simplifies to
\be
\CM_P(\{C_r\},\ul \Lambda;\ul x)=M_P(\{d_r^*\}, \Lambda; x),
\ee
with 
\be
\label{eq:DefMP}
\begin{split}
M_P(\{c_r\},\Lambda;x)&=\sqrt{|\Delta(\{c_j^\star\})|}\left(\frac{i}{\pi }\right)^{P}\\
&\qquad \times \int_{\langle
 \{c_r\}\rangle-ix^{||}}\prod_{j=1}^{P} \frac{1}{ B(c_j^\star,z)}
e^{-\pi  Q(z)-2\pi i  B(z,x)} d^{P}z.
\end{split}
\ee
 
We continue with the part in Eq. (\ref{eq:PhiHatML}) that depends on
$x_{s_j}$ and demonstrate that the sum over these integers can be
carried out as a geometric sum. The kernel contains terms of the form 
\be
\sgn(\underline B(C_{s}^{\perp V},\underline k))=\sgn\!\left(x_s-\sum_{j=1}^L
c_{s,v_j} x_{v_j}\right),
\ee 
which we can write equivalently as
\be
\sgn\!\left(\underline B(C_{s},\underline k+(0,\nu^{||})\right)=\sgn\!\left(\underline
B\left(C_{s},(k,\sum_{j=1}^L x_{s_j} d_{s_j}+\nu^{||})\right)\right),
\ee 
with $(k,\sum_j x_{s_j} d_{s_j}+\nu^{||})\in \Lambda^*\oplus
(-\Lambda_d)^*$, and $\nu^{||}$
\be
\label{eq:nu||}
\nu^{||}=-\sum_{j=1}^{M-L}\sum_{\ell=1}^L c_{s_j,v_\ell} x_{v_\ell} d_{s_j}.
\ee
We can now rewrite the term $B(x_rd_r,k)$ in the exponent of $q$ as
\be
B(x_rd_r,k)=B\!\left(\sum_{\ell=1}^L x_{v_\ell} d_{v_\ell} -\nu^{||}+\sum_{j=1}^{M-L} x_{s_j} d_{s_j}+ \nu^{||},k\right).
\ee  
We prove in Appendix \ref{sec:ProofId} the orthogonality relation
\be
\label{eq:orthorelation}
B\!\left(d_{s_j},d_{v_\ell}+\sum_{k=1}^{M-L}
c_{s_k,v_\ell}d_{s_k}\right)=0.
\ee
Thus the vector $d_{v_\ell}+\sum_k
c_{s_k,v_\ell}d_{s_k}$ is the component of $d_{v_\ell}$ orthogonal
to the lattice $\Lambda(\{d_{s_j}\})$,
\be 
d_{v_{\ell}}^{\perp S}=d_{v_\ell}+\sum_k
c_{s_k,v_\ell}d_{s_k}.
\ee 
Consider the substitution
\be
\label{eq:ktrafo}
k = k'-\sum_{\ell=1}^L x_{v_\ell}(d_{v_\ell}+\sum_j
c_{s_j,v_\ell}d_{s_j})\in \Lambda+\mu,
\ee
such that
\be
\label{eq:kprime}
k'\in \Lambda+\mu-\nu^{||}.
\ee
This substitution expresses the exponent of $q$ in Eq. (\ref{eq:PhiHatML}) as
\be
Q(k)/2+B(x_rd_r,k)= -Q(x_{v_\ell} d_{v_{\ell}}^{\perp S})/2 + Q(k')/2+B(\sum_{j} x_{s_j} d_{s_j}+ \nu^{||},k').
\ee
Moreover, this substitution changes the terms with the elliptic variables $u$ and $v$ to
\be
B(v,k)+ B(x_rd_r,u) = 
  B(u-v,x_{v_\ell} d_{v_\ell}^{\perp S})+ B(v,k')+ B(u,\nu^{||}+x_{s_j}d_{s_j}).
\ee 

Comparing with the definition of $\Phi^+_{\mu,\nu}$
(\ref{eq:DefPhi+}), we realize that the sum over $x_{s_j}$ combines to
$\Phi^+_{\mu,\nu^{||}}(\tau,u^{||},v,\{d_{s_j}\})$, with $\nu^{||}$
and $u^{||}$ the components of $\nu$ and $u$ parallel to the hyperplane
spanned by $\{d_{s_j}\}$.
We thus arrive at 
\be
\label{eq:PhiHatML2} 
\begin{split} 
&\CR_\mu(\tau,\bar \tau, u,\bar u, v, \bar v,\{d_r\})=
 \sum_{L=1}^{M} \sum_{\{v_1,\dots,v_L,s_{1},\dots,
  s_{M-L}\}\in \{ 1,\dots,M\}}\\
&\qquad    \sum_{x_{v_\ell}\in \mathbb{Z} } 2^{-L} M_L(\{d_{v_i}^*\}, \Lambda; 
  \sqrt{2y}( x_{v_\ell} +{\rm Im}(\sigma_{v_\ell})/y) d_{v_\ell})\, \\
&\qquad \times q^{-Q(x_{v_\ell} d_{v_{\ell}}^{\perp S})/2}\, e^{2\pi i
  B(u-v,x_{v_\ell} d_{v_\ell}^{\perp S})} \\
&\qquad \times \Phi^+_{\mu,\nu^{||}}(\tau,u^{||},v,\{d_{s_j}\}).
\end{split}
\ee
Thus $\widehat \Phi_\mu=\Phi^+_\mu+\CR_\mu$ is fully expressed in terms of
$\Phi^+_{\mu,\nu^{||}}$ for different sets $\{d_{s_j}\}$ and data in terms of $\Lambda$.

\subsubsection*{The Lattice $\Lambda_d(\{d^*_{v_\ell}\})$ and Glue Vectors}
We next decompose the sum over $x_{v_\ell}$ into a sum over an
integral lattice together with a finite sum over conjugacy classes.
The vectors $\sum_\ell x_{v_\ell}
d_{v_\ell}^{\perp S}$ lie in the vector space over $\mathbb{R}$ orthogonal to
$\{d_{s_j}\}$. While this space is spanned by $d_{v_\ell}^{\perp S}$,
$\ell=1,\dots,L$, it is also spanned by the dual basis vectors
$d_{v_\ell}^*$ since $B(d_{v_\ell}^*,d_{s_j})=0$ for all $\ell=1,\dots,
L$ and $j=1,\dots,M-L$. We define $\Lambda_d(\{d^*_{v_\ell}\})\subset
\Lambda_d$ as the {\it integral}
sublattice of $\Lambda_d$ spanned by $d_{v_\ell}^*$. The
generators of $\Lambda_d(\{d^*_{v_\ell}\})$ are thus suitable linear
combinations of $d^*_{v_\ell}\in
\Lambda_d^*$ such that $\sum_{\ell} N_\ell d^*_{v_\ell}\in \Lambda_d$. 

Because of the orthogonality of the sets
$\{d^*_{v_\ell}\}$ and $d_{s_j}$, we have the equivalence
$\Lambda_d(\{d^*_{v_\ell}\})\equiv \Lambda_d(\{d^{\perp
  S}_{v_\ell}\})$. Moreover the direct sum of the integral lattices is a
sublatice of $\Lambda_d$, $\Lambda_d(\{d_{s_j}\})\oplus \Lambda_d(\{d^*_{v_\ell}\})\subset
\Lambda_d$, while $\Lambda_d$ is a sublattice of the direct sum of the
dual lattices $\Lambda_d\subset \Lambda_d^*(\{d_{s_j}\})\oplus \Lambda^*_d(\{d^*_{v_\ell}\})$. Using general theory of lattices \cite{conway}, an element 
$k\in \Lambda_d$ can be written as
$k=\nu_g+\ell^{||}+\ell^\perp$ with $\ell^{||}\in
\Lambda_d(\{d_{s_j}\})$ and $\ell^\perp\in \Lambda_d(\{d^*_{v_\ell}\})$ and $\nu_g\in
\Lambda_d$ a glue vector.
The number $\mathcal{N}_g(\{d_{s_j}\},\Lambda_d)$ of glue vectors $\nu_g$, $g=1,\dots, \mathcal{N}_g$, is given in terms of the number of
elements in the cosets of the lattices involved,
\be
\label{eq:CNg}
\mathcal{N}_g(\{d_{s_j}\},\Lambda_d)=\sqrt{\frac{|\Lambda^*_d(\{d^*_{v_i}\})/\Lambda_d(\{d^*_{v_i}\})|\, |\Lambda^*_d(\{d_{s_j}\})/\Lambda_d(\{d_{s_j}\})|}{|\Lambda_d^*/\Lambda_d|}}.
\ee
We have the orthogonal decomposition $\nu_g=\nu_g^{||}+\nu_g^{\perp}$, with
$\nu^{||}_g$ the projection of $\nu_g$ to $\Lambda^*_d(\{d_{s_j}\})$
and $\nu_g^{\perp}$ the projection to $\Lambda_d(\{d^*_{v_\ell}\})$. The number $\mathcal{N}_g$ is easily
determined by using that for a general lattice
$\Lambda$, the number of elements $|\Lambda^*/\Lambda|$ is determined
in terms of its bilinear form $B$ as, 
\be
|\Lambda^*/\Lambda|=|{\rm det}(B)|.
\ee
 
Now we return to the completion \eqref{eq:PhiHatML2}. Recall that
$M_L(\{d_{v_\ell}^*\}, \Lambda;k)$ only depends on the components of
$k$ parallel to $\{d_{v_\ell}^*\}$. The vectors
$x_{v_\ell}d_{v_{\ell}}$ on the second line
of \eqref{eq:PhiHatML2} can thus be replaced by
$x_{v_\ell}d_{v_{\ell}}^{\perp S}$ as on the third line. Moreover, since  $\{d_{v_\ell}^*\}$ generates
the same lattice as $d_{v_\ell}^{\perp S}$, the sum over $x_{v_\ell}$
can be written as a sum over $\Lambda_d(\{d_{v_\ell}^*\})+\nu$ for
specific $\nu$.

\subsubsection*{Structural Formula for $\widehat \Phi_\mu$}
We combine the above results, to arrive at a structural formula for
$\widehat \Phi_\mu$. To this end, let $\{d_r\}$ be a set with $M\geq
0$ linearly independent elements $d_r\in \Lambda$ spanning the lattice
$\Lambda_d\subset \Lambda$, and with dual vectors $d_r^* \in
\Lambda_d^*$. Let $\{d_{s_j}\}\subsetneq \{d_r\}$ be subsets
with $M-L<M$ elements with $j=1,\dots, M-L$. Moreover, let $\nu_g$,
$g=1,\dots, \mathcal{N}_g$, be the glue vectors for gluing of the orthogonal sublattices $\Lambda(\{d_{s_j}\})$ and
$\Lambda(\{d_r^*\}/\{d^*_{s_j}\})$ in $\Lambda_d$ as discussed
above.  

The modular completion
$\widehat \Phi_{\mu}$ of $\Phi_\mu$ then reads,
\be 
\label{eq:StructForm}
\boxed{
\begin{split}
& \widehat 
 \Phi_\mu(\tau,\bar \tau,u,\bar u,v,\bar v,\{d_r\})= \Phi^+_\mu(\tau,u,v,\{d_r\})\\
 &\qquad +
\sum_{L=1}^M \sum_{\{d_{s_j}\}\subsetneq \{d_r\}} \sum_{g=1}^{\mathcal{N}_g} 2^{-L}
 R_{L,\nu_g^\perp}(\{d^*_r\}/\{d^*_{s_j}\},\Lambda; \tau,\bar
 \tau,u^\perp-v^\perp,\bar u^\perp-\bar v^\perp) \\
& \qquad \quad \times \Phi^+_{\mu,\nu_{g}^{||}}(\tau,u^{||},v,\{d_{s_j}\}),
\end{split}}
\ee
where the non-holomorphic function $R_{L,\nu}$ is defined by 
\be
\begin{split}  
R_{L,\nu}(\{d^*_v\},\Lambda ;\tau,\bar \tau,\sigma,\bar \sigma)&=\sum_{k\in
  \Lambda_d(\{d^*_{v}\})+\nu}
M_L( \{d^*_{v}\}, \Lambda; \sqrt{2y}(k+{\rm Im}(\sigma)/y ))\,\\
&\qquad \times q^{-Q(k)/2}e^{2\pi i B(\sigma, k)},
\end{split}
\ee 
with $M_L$ as in Eq. (\ref{eq:DefMP}). Using the relation between
$\Phi_\mu=\Phi_{\mu,0}$ and $\Phi_{\mu,\nu}$ (\ref{eq:PhimunuPhimu0})
and the relation (\ref{eq:Phimum}), we determine from Eq. (\ref{eq:StructForm}) that the completion
$\widehat \Phi_{\mu,\nu}$ of  $\Phi_{\mu,\nu}$ is given by 
\be 
\label{eq:StructForm2}
\boxed{
\begin{split}
& \widehat 
 \Phi_{\mu,\nu}(\tau,\bar \tau,u,\bar u,v,\bar v,\{d_r\})= \Phi^+_{\mu,\nu}(\tau,u,v,\{d_r\})\\
 & +
\sum_{L=1}^M \sum_{\{d_{s_j}\}\subsetneq \{d_r\}} \sum_{g=1}^{\mathcal{N}_g} 2^{-L}
 R_{L,\nu_g^\perp+\nu^\perp}(\{d^*_r\}/\{d^*_{s_j}\},\Lambda; \tau,\bar
 \tau,u^\perp-v^\perp,\bar u^\perp-\bar v^\perp) \\
& \quad \times \Phi^+_{\mu,\nu_{g}^{||}+\nu^{||}}(\tau,u^{||},v,\{d_{s_j}\}).
\end{split}}
\ee

The functions $M_L( \{d^*_{v}\}, \Lambda; \sqrt{2y}(k+{\rm
  Im}(\sigma)/y ))$ vanish in the limit $y\to \infty$, except when the last
argument vanishes, $k+{\rm Im}(\sigma)/y=0$. In the latter case, they actually contribute
a holomorphic term. For example for $L=2$,
$M_2(\{c_j\},\Lambda;0)=(2/\pi)\,{\rm Arctan}(\alpha)$ with
$\alpha=B(c_1,c_2)/\sqrt{\Delta(c_1,c_2)}$ \cite[Eq. (3.23)]{Alexandrov:2016enp}. 
Finally, we can make use of the periodicity \eqref{eq:Phi+period}, and
Eq. \eqref{eq:PhiPlusPhi} to express $\widehat \Phi$ in terms of the
original Appell functions $\Phi$.

\section{Application to Appell Functions for BPS Indices}
\label{sec:AppellBPS}

We give in this section various examples of the general results in
previous sections based on the building blocks
$\Psi_{(r_1,\dots,r_\ell),(a,b)}$ (\ref{psi}) for generating
functions of certain BPS indices.
They serve to illustrate as well as to verify the general
discussion. 

\subsection{The Lattice $A_2$} 
This case has been studied in some detail in
Ref. \cite{Alexandrov:2016enp, Manschot:2017xcr}. To make the connection, we make the
change $k_1\to -k_1$ in \cite[Eq. (5.21)]{Manschot:2017xcr}.
The Appell-Lerch sums defined for the $A_2$ lattice then have
\be
d_1=\left(\begin{array}{c} -1 \\ 0 \end{array}\right),\qquad d_2=\left(\begin{array}{c} 1 \\ 1 \end{array}\right),\qquad B(d_1,d_2)=-1.
\ee
The dual vectors are
\be
d_1^*=\frac{1}{3}\left(\begin{array}{c} -1 \\
                         1 \end{array}\right),\qquad
                     d_2^*=\frac{1}{3}\left(\begin{array}{c} 1 \\
                                              2 \end{array}\right),\qquad
                                          Q(d_1^*)=Q(d_2^*)=\frac{2}{3}.
\ee 
Moreover the orthogonal projections are
\be
d_1^{\perp 2}= \frac{1}{2}\left(\begin{array}{c} -1 \\
                         1 \end{array}\right),\qquad d_2^{\perp
                       1}=\frac{1}{2}\left(\begin{array}{c} 1 \\
                                             2 \end{array}\right),\qquad 
Q(d_1^{\perp 2})=                 
Q(d_2^{\perp 1})=\frac{3}{2}.
\ee 
These are proportional to the dual vectors $d_r^*$, such that $d_r^*$
and $d_r^{\perp s}$ generate the same integral lattice. Moreover,
$d_1$ and $d_2$ generate the lattice $\Lambda$, such that 
$\Lambda(\{d_1,d_2\})=\Lambda_d\equiv \Lambda$. 

The extended quadratic form ${\underline {\bf A}}$ reads
\be
{\underline {\bf A}}=\left( \begin{array}{cccc} 2 & -1 & 2 & -1 \\ -1 &
                                                                       2
                                                       & -1 & -1 \\ 2
                                                  & -1 & 0 & 0 \\ -1 & -1 & 0 & 0 \\  \end{array} \right).
\ee 
The vectors ${\bf C}_r$ and ${\bf C}_r'$ are 
\be
{\bf C}_1= -\frac{1}{3}\left(\begin{array}{c} -1 \\ 1 \\ 2 \\ 1 \end{array}\right) ,\quad {\bf C}_1'=
\left(\begin{array}{c} 0 \\ 0 \\ 1 \\ 0 \end{array}\right) ,\quad {\bf
  C}_2=-\frac{1}{3}
\left(\begin{array}{c} 1\\ 2\\ 1\\ 2 \end{array}\right) ,\quad {\bf C}_2'=
\left(\begin{array}{c} 0 \\ 0 \\ 0 \\ 1  \end{array}\right).
\ee  
Note ${\bf C}_r\in \underline \Lambda^*$, it is not an element of
$\underline \Lambda$. We have for $c_{1,2}$ \eqref{eq:csv}
\be 
c_{1,2}=\frac{1}{2}.
\ee
Let us consider the term corresponding to $L=1$ in Eq. (\ref{eq:PhiHatML2}) with $v_1=1$ and
$s_1=2$. We verify the general relation \eqref{eq:orthorelation} that $d_2+c_{12} d_1$ is perpendicular to $d_1$.
 
The integral lattice $\Lambda(d_1^*)$ is generated by $3d_1^*\in \Lambda$, and
$|\Lambda^*(d_1^*)/\Lambda(d_1^*)|=6$, such that the elements of
$\Lambda^*(d_1^*)/\Lambda(d_1^*)$ are $j/2\times
d_1^*$ with $j=0,\dots,5$. Moreover, since
$|\Lambda^*(d_2)/\Lambda(d_2)|=2$ and $|\Lambda^*/\Lambda|=3$, we
determine with Eq. \eqref{eq:CNg} that there are $\mathcal{N}_g=2$
glue vectors. The glue vectors are $\nu_1=0$ and  $\nu_1=d_1\in
\Lambda$ with $\nu_1^{||}\in \Lambda^*(d_2)$ and $\nu_1^{\perp}\in
\Lambda^*(d_1^*)$.

With these data we can indeed verify that the modular completions of
these Appell functions derived in \cite[Thm 5.3 and Eq. (5.27)]{Alexandrov:2016enp} and \cite[Eq. (5.24)
and (5.25)]{Manschot:2017xcr} indeed have the form of Eq. (\ref{eq:StructForm}).

\subsection{The Lattice $A_3$} 
\label{sec:LatA3}
We study the case of the Appell function constructed from the $A_3$ lattice 
as in Eq. (\ref{psi}),
\begin{eqnarray}
\Psi_{(1,1,1,1),(0,0)}(\tau,z) &=& \sum_{b_1,b_2,b_3\in\mathbb{Z}}\frac{w^{-6b_1-4b_2-2b_3}q^{b_1^2+b_2^2+b_3^2+b_1b_2+b_2b_3+b_1b_3}}{(1-w^4 q^{b_1-b_2})(1-w^4 q^{b_2-b_3})(1-w^4 q^{b_1+b_2+2b_3})},
\end{eqnarray}
with $w=e^{2\pi i z}$ as before. To bring the quadratic form to the
standard form for the lattice $A_3$ with the simple roots as the basis vectors, we make the
following transformation,
\begin{eqnarray}
b_1\to k_2-k_1, \quad b_2\to k_1,\quad b_3\to k_3-k_2.
\end{eqnarray}
This brings the quadratic form to the desired form for $A_3$,
\begin{eqnarray}
\begin{pmatrix}
2 & 1 & 1\\
1 & 2 & 1\\
1 & 1 & 2
\end{pmatrix}\rightarrow \begin{pmatrix}
2 & -1 & 0\\
-1 & 2 & -1\\
0 & -1 & 2
\end{pmatrix},
\end{eqnarray}
such that $\Psi_{(1,1,1,1),(0,0)}$ reads after the change of
summation variables,
\be
\label{eq:PsiA3}
\begin{split}
&\Psi_{(1,1,1,1),(0,0)}(\tau,z)=\\
&\qquad\sum_{k_1,k_2,k_3\in\mathbb{Z}^3}\frac{w^{2k_1-4k_2-2k_3}q^{k_1^2+k_2^2+k_3^2-k_1k_2-k_2k_3}}{(1-w^4
  q^{-2k_1+k_2})(1-w^4 q^{k_1+k_2-k_3})(1-w^4 q^{-k_2+2k_3})}.
\end{split}
\ee
This is of the form $\Phi_{\mu,\nu}(\tau,u,v,\{d_r\})$. We note that the choice $a=b=0$ in $\Psi_{(r_1,\dots,r_\ell),(a,b)}$
(\ref{psi}) corresponds to $\mu=\nu=0$ in $\Phi_{\mu,\nu}$. In the
remainder of this section, we will abbreviate
$\Psi_{(1,1,1,1),(0,0)}=:\Psi$. 

We can now read of the vectors $d_j$,
\be
d_1=\left(\begin{array}{c} -1 \\ 0 \\ 0 \end{array} \right),  \quad
                                        d_2=\left(\begin{array}{c}
                                                        1 \\ 1 \\ 0 \end{array}
                                                                 \right),
                                                                \quad d_3=\left(\begin{array}{c}
                                                        0 \\ 0 \\ 1\end{array}
                                                                 \right).
                                                                 \ee
These vectors are related to the basis of positive roots
$\{\alpha_i\}$ by the Weyl reflection $S_{\alpha_1}$
\eqref{eq:Salpha}. This reflection equals here ${\bf D}^{-1}{\bf C}$,
with $\bf C$ defined below Eq. \eqref{eq:ulA} and $\bf D$ defined below
Eq. \eqref{DefB}.
                                                                 
The matrix of innerproducts of these vectors is identical to the
$A_{3}$ root lattice as we can write $d_j=S_{\alpha_1}(\alpha_j)$.                                                               
The dual vectors are
\be
d_1^*=\left(\begin{array}{r} -\frac{1}{4} \\ \frac{1}{2} \\ \frac{1}{4} \end{array} \right),  \quad
                                        d_2^*=\left(\begin{array}{c}
                                                        \frac{1}{2} \\ 1 \\ \frac{1}{2} \end{array}
                                                                 \right),
                                                                \quad d_3^*=\left(\begin{array}{c}
                                                        \frac{1}{4} \\ \frac{1}{2} \\ \frac{3}{4}\end{array}
                                                                 \right).
\ee  

Comparison of Eq. \eqref{phigen2} and Eq. \eqref{eq:PsiA3} demonstrates that the elliptic
variables $u$ and $v$ are written in terms of the basis as
\be 
u=4z\,(d_1^*+d_2^*+d_3^*)=z \left(\begin{array}{c} 2 \\ 8 \\ 6 \end{array}\right),\qquad v=-z \left(\begin{array}{c} 1 \\ 4 \\ 3 \end{array}\right).
\ee  

Expanding the geometric sums, we obtain for $\Psi$, 
\be
\begin{split}
\Psi(\tau,z) &= \sum_{k_j
                                   \in\mathbb{Z}} K^\epsilon(\{k_j\}, a)\, w^{2k_1-4k_2-2k_3+4(k_4+k_5+k_6)},\\
&\,\times
q^{k_1^2+k_2^2+k_3^2-k_1k_2-k_2k_3+k_4(k_2-2k_1)+k_5(k_1+k_2-k_3)+k_6(2k_3-k_2)},
\end{split}
\ee
with the kernel $K^\epsilon(\{k_j\}, a)$ given by 
\be
\begin{split} 
K^\epsilon(\{k_j\}, a)&=\tfrac{1}{8}({\rm sgn} (k_2-2k_1+4a)+{\rm sgn}(k_4+\epsilon))\\ 
                                 &\,\times ({\rm
                                   sgn}(k_1+k_2-k_3+4a)+{\rm
                                   sgn}(k_5+\epsilon))\\
                                 &\,\times ({\rm sgn}(2k_3-k_2+4a)+{\rm sgn}(k_6+\epsilon)),
\end{split}
\ee
and $a={\rm Im}(z)/y$.\footnote{Note $a$ here is different from the
  use of $a$ in Eq. \eqref{psi}.} The six-dimensional sum $k_j\in \mathbb{Z}$
is identified with the sum over $\ul k \in \ul \Lambda$. Either by direct
computation or using Eqs
(\ref{eq:zrhosigma}),  \eqref{eq:sigma} and \eqref{eq:rho}, we deduce for the
elliptic variable $\ul z\in \ul \Lambda \otimes \mathbb{C}$,
\be
\label{eq:uzA3}
\ul z=z\left(\begin{array}{c} 2 \\ 8 \\ 6 \\ -9 \\ -12 \\ -9  \end{array}\right)_{\alpha\gamma}.
\ee 
Comparison with Eq. \eqref{eq:zrhosigma} demonstrates that
$\rho=u=z\,(2,8,6)^T$, while $\sigma=-z\,(9,12,9)^T$.

We express the completion $\widehat \Psi$ of $\Psi$ as the sum
\be
\widehat \Psi(\tau,\bar \tau,z,\bar z)=\Psi(\tau,z)+\mathcal{R}(\tau,\bar \tau,z,\bar z).
\ee
Moreover, we let $\Psi^+$ be the function with the $\epsilon$ in the kernel replaced
by the appropriate shift of $a$ as for $\Phi^+_{\mu,\nu}$
(\ref{eq:DefPhi+}), and $\mathcal{R}$ the non-holomorphic term to be determined. 

We now introduce alternative notation for the arguments of $E_P$ and $M_P$
using equivalent analytic expressions.  We recall from \cite{Alexandrov:2016enp}
\be 
E_1(\{c\},\Lambda;x)\to E_1(u_1)\quad {\rm with}\quad u_1= \frac{B(c,x)}{\sqrt{Q(c)}},
\ee
and equal to the error function. For $E_2$, we introduce
\be
E_2(\{c_1,c_2\},\Lambda;x)\to E_2(\alpha;u_1,u_2),
\ee 
with 
\be
\alpha=\frac{B(c_1,c_2)}{\sqrt{\Delta(c_1,c_2)}},\quad
u_1=\frac{B(c_{1\perp 2},x)}{\sqrt{Q(c_{1\perp 2})}},\quad u_2=\frac{B(c_{2},x)}{\sqrt{Q(c_{2})}},
\ee
and $E_2$ given by \cite[Eq. (3.29)]{Alexandrov:2016enp}.

We can now write $\mathcal{R}$ in terms of the generalized error
functions $E_1,E_2,E_3$ using the value of $\ul z$,
{\scriptsize{\be\label{eq:CRA4}\begin{split} 
& \mathcal{R}(\tau,\bar \tau,z, \bar z )=\\
 &\frac{1}{8} \sum_{k_j\in \mathbb{Z}} \left[\left({\rm sgn}(k_4-9a)-E_1(\frac{2\sqrt{2y}}{\sqrt 3}(k_4-9a))\right){\rm sgn}(k_1+k_2-k_3+4a){\rm sgn}(2k_3-k_2+4a)\right.\\ 
&+\left. \left({\rm sgn}(k_6-9a)-E_1(\frac{2\sqrt{2y}}{\sqrt 3}(k_6-9a))\right){\rm sgn}(k_1+k_2-k_3+4a){\rm sgn}(-2k_1+k_2+4a)\right.\\ 
&+\left. \left({\rm sgn}(k_5-12a)-E_1(\sqrt{2y}(k_5-12a))\right){\rm sgn}(-2k_1+k_2+4a){\rm sgn}(2k_3-k_2+4a)\right.\\ 
&+ \left.{\rm sgn}(2k_3-k_2+4a)\left({\rm sgn}(k_4-9a){\rm sgn}(k_5-12a)-E_2\left(\frac{1}{\sqrt{2}},2\sqrt{y}(k_4-k_5/2-3a),\sqrt{2y}(k_5-12a)\right)\right) \right.\\ 
&+\left. {\rm sgn}(-2k_1+k_2+4a)\left({\rm sgn}(k_6-9a){\rm sgn}(k_5-12a)-E_2\left(\frac{1}{\sqrt{2}},2\sqrt{y}(k_6-k_5/2-3a),\sqrt{2y}(k_5-12a)\right)\right)\right.\\ \
&+ \left. {\rm sgn}(k_1+k_2-k_3+4a)\left({\rm sgn}(k_4-9a){\rm sgn}(k_6-9a)-E_2\left(\frac{1}{\sqrt{8}},\sqrt{3y}(k_4-k_6/3-6a),\frac{2\sqrt{2y}}{\sqrt 3}(k_6-9a)\right)\right)\right.\\ 
&+\left. {\rm sgn}(k_4-9a){\rm sgn}(k_5-12a){\rm
    sgn}(k_6-9a)-\mathcal{E}_3(\{C_1,C_2,C_3\},\ul \Lambda; \sqrt{2y}(\ul k +\ul a))\right]\\
&\times     q^{\ul{Q}(\ul{k})/2}e^{2\pi
    i \ul{B}(\ul{k},\ul z)}.
\end{split}\ee}}
where
\be
\label{eq:aula}
a=\frac{{\rm Im}(z)}{y},\quad {\rm and} \quad \ul a=\frac{{\rm Im}(\ul z)}{y}.
\ee
We write $\mathcal{R}$ as a sum of three terms involving either $M_1$,
$M_2$ or $M_3$,
\be
\label{eq:CRK123}
\mathcal{R}=\sum_{\ul k\in \ul \Lambda}
\left( \widehat K_1+ \widehat K_2+ \widehat M_3(\{C_1,C_2,C_3\},\ul k+\ul a)\right)\,q^{\ul{Q}(\ul{k})/2} e^{2\pi i
  \ul{B}(\ul{k}, \ul z)},
\ee  
and discuss each of these terms separately. To express
$\Phi^+_{\mu,\nu}$ as an Appell function $\Phi_{\mu,\nu}$, we will
choose the imaginary part of the elliptic variable negative and
sufficiently small, such that
\be
\label{eq:-1a0}
-1\ll a<0,
\ee
such that the components $\nu_r$ and $\sigma_r$ satisfy
\be
\lfloor \nu_r + {\rm Im}(\sigma_r)/y \rfloor=\lfloor \nu_r \rfloor,
\ee
for all $\nu=\sum_r \nu_r d_r \in \Lambda_d^*$. As a result, these
terms in Eq. (\ref{eq:Phi+Appell}) simplify for this specific case to
the fractional part $\nu_r-\lfloor \nu_r \rfloor=\{\nu_r\}$ of the
components $\nu_r$. We introduce $\tilde \nu$,
\be
\tilde \nu= \sum_r \{\nu_r\} d_r,
\ee
such that for this choice of $z$,
\be
\Phi^+_{\mu,\nu}=\Phi_{\mu,\tilde \nu}.
\ee
Recall $\Phi^+_{\mu,\nu}$ is periodic in $\nu$ under shifts by an element in
$\Lambda_d$, while $\Phi_{\mu,\tilde \nu}$ is not.

We discuss below in detail the various terms of the completion with
$L=1,2,3$. These terms contribute to the full non-holomorphic
completion of the $SU(4)$ partition function of VW-theory
\cite{Manschot:2014cca}. Indeed, the terms below reproduce various terms given in
\cite[Appendix F.3]{Alexandrov:2019rth}, which conjectured the
completion of the $SU(4)$ partition function from a
string theory perspective. Moreover, the non-holomorphic term
$\mathcal{R}$ (\ref{eq:CRA4}) is in agreement with the structure of
the non-holomorphic part of refined partition functions derived in
\cite[Theorem 1]{Alexandrov:2020bwg}. We leave further analysis for
future work.

\subsection*{Terms with $L=1$}
We first consider the term $\widehat K_1$, which is the sum of terms
involving a single $M_1$. We have
\be
\label{eq:PsiK1}
\begin{split}
\widehat K_1&=\tfrac{1}{2}\kappa_1\, M_1(2\sqrt{2y/3}\,(k_4-9a))+\tfrac{1}{2}\kappa_2\,
M_1(\sqrt{2y}(k_5-12a))\\
&\quad +\tfrac{1}{2}\kappa_3\,M_1(2\sqrt{2y/3}(k_6-9a)),
\end{split}
\ee
where, 
\be
\label{eq:kappas}
\begin{split}
  \kappa_1 &= \tfrac{1}{4} (\sgn(2k_3-k_2+4a)+\sgn(k_6-k_4/3-6a))\\
  &\quad \times (\sgn(k_1+k_2-k_3+4a)+\sgn(k_5-2k_4/3-6a)),\\ 
  \kappa_2 &= \tfrac{1}{4} (\sgn(-2k_1+k_2+4a)+\sgn(k_4-k_5/2-3a))\\
&\quad \times (\sgn(2k_3-k_2+4a)+\sgn(k_6-k_5/2-3a)),\\ 
\kappa_3 &= \tfrac{1}{4} (\sgn(-2k_1+k_2+4a)+\sgn(k_4-k_6/3-6a))\\
&\quad \times (\sgn(k_1+k_2-k_3+4a)+\sgn(k_5-2k_6/3-6a)).
\end{split}
\ee

To compare with the general discussion of Section \ref{sec:ModCom}, we
determine the coefficients $c_{i,j}$ \eqref{eq:csv}, 
\be
\begin{split}
c_{2,1}=\frac{2}{3},\qquad c_{3,1}=\frac{1}{3},\\
c_{1,2}=\frac{1}{2},\qquad c_{3,2}=\frac{1}{2},\\
c_{1,3}=\frac{1}{3},\qquad c_{2,3}=\frac{2}{3}.\\
\end{split}
\ee

We now address each of the terms in Eq. (\ref{eq:PsiK1}), 
\begin{enumerate}         
\item To carry out the sum over $k_5,k_6\in \mathbb{Z}$ in the term with $\kappa_1$, we first replace
$k_5= k'_5+2k_4/3, k_6= k'_6+k_4/3$. In the
general discussion of Section \ref{sec:ModCom}, this is the
shift of the vectors $\sum_j x_{s_j}d_{s_j}$ by $\nu^{||}$ \eqref{eq:nu||}.
The next step is the following substitution 
\be k_1= k'_1+k_4/3, \qquad k_2= k'_2-2k_4/3, \qquad
k_3= k'_3-k_4/3,\ee  
which corresponds to the substitution \eqref{eq:ktrafo} in Section
\ref{sec:ModCom}. Since $k'\in \Lambda+\mu-\nu$ (\ref{eq:kprime}), we deduce that
for $k_4\in 3\mathbb{Z}+g, g=0,1,2,$
\be 
\nu^{||}_g=-\frac{g}{3}\left(\begin{array}{c} -1 \\ 2 \\
                          1 \end{array}\right) \mod \mathbb{Z}.
\ee

The lattice $\Lambda(d_1^*)$ is generated by
$4d_1^*$. The number of glue vectors is $\mathcal{N}_g=\sqrt{12\times
  3/4}=3$. These can be choses as $\nu_0=0, \nu_1=d_1$ and $\nu_2=2d_1$ matching the sum over
$g=0,1,2$. We get $\Phi^+_{0,\nu_g^{||}}(\tau, u,v,\{d_2,d_3\})$,
\be
\begin{split}
&\sum_{k_j\in \mathbb{Z}}\kappa_1\, M_1(2\sqrt{2y/3}(k_4-9a))\, q^{\ul
  Q(\ul k)}e^{2\pi i\ul B(\ul k,\ul z)} =\\
& \times \sum_{k_4\in 3\mathbb{Z}+g \atop g=0,1,2}
                                                   q^{-2k_4^2/3}w^{12k_4}\,M_1(2\sqrt{2y/3}(k_4-9a))
                                           \,\Phi^+_{0,\nu^{||}_g}(\tau, u,v,\{d_2,d_3\}).
\end{split}
\ee
We evaluate the sums over $k_5$ and $k_6$ for a specific choice of $a$ and $\nu_g$.
Namely with $a$ as in Eq. (\ref{eq:-1a0}), and
\be
\tilde \nu^{||}_0=0,\qquad 
\tilde \nu^{||}_1= \frac{1}{3}\left(\begin{array}{c} 1 \\ 1 \\
                          2 \end{array}\right), \qquad 
                      \tilde \nu^{||}_2= \frac{1}{3}\left(\begin{array}{c} 2 \\ 2 \\
                          1 \end{array}\right),
                      \ee
one can show that $S_{0,\tilde \nu^{||}_g}$ vanishes. Then $\Phi^+_{0,\nu^{||}_g}$ agrees with the Appell function
$\Phi_{0, \tilde \nu^{||}_g}$ \eqref{phigen2}.

\item To carry out the sum over $k_4,k_6\in \mathbb{Z}$ in the term with $\kappa_2$, we first make the
  substitution $k_4= k'_4+k_5/2$, and $k_6= k_6'+k_5/2$
  \eqref{eq:nu||}. Furthermore, we make the substitution \eqref{eq:ktrafo}
  \be
  k_1= k'_1-k_5/2, \qquad k_2= k'_2-k_5,\qquad k_3=k'_3-k_5/2,
  \ee
  such that for $k_5\in 2\mathbb{Z}+g, g=0,1$,
\be 
\nu^{||}_g=-\frac{g}{2}\left(\begin{array}{c} 1 \\ 2 \\
                          1 \end{array}\right) \mod \mathbb{Z}.
\ee

  The lattice $\Lambda(d_2^*)$ is generated by $2d_2^*$.
The number of glue vectors is in this case, $\mathcal{N}_g=\sqrt{4\times
  4/4}=2$, for which we take 0 and $d_2$. The sum then evaluates to
\be
\begin{split}
&\sum_{k_j\in \mathbb{Z}}\kappa_2\, M_1(\sqrt{2y}(k_5-12a))\, q^{\ul
  Q(\ul k)}e^{2\pi i\ul B(\ul k,\ul z)} =\\
& \times \sum_{k_5\in 2\mathbb{Z}+g \atop g=0,1}
                                                   q^{-k_5^2/2}w^{12k_5}\,M_1(2\sqrt{2y}(k_5-12a))
                                           \,\Phi^+_{0,\nu^{||}_g}(\tau,
                                           u,v,\{d_1,d_3\}).
                                         \end{split}
                                         \ee
                                         As before, we relate
                                         $\Phi^+_{\nu^{||}_g}$ to an Appell
                                         function. With $a$ as in Eq. (\ref{eq:-1a0}), and
                                         \be
                                         \tilde \nu^{||}_g=\frac{g}{2}\left(\begin{array}{c} -1 \\ 0 \\
                          1 \end{array}\right),\quad g=0,1,
                      \ee
                      then $\Phi^+_{\nu^{||}_g}$ agrees with
                      $\Phi_{\tilde \nu^{||}_g}$.

\item Finally for $L=1$, to carry out the sum over $k_4,k_5\in \mathbb{Z}$ in
  the term with $\kappa_3$, we make the substitution
  $k_4= k'_4+k_6/3$, $k_5= k'_5+2k_6/3$ \eqref{eq:nu||}, followed by
  the substitutions
\be
  k_1= k_1'-k_6/3, \quad k_2= k_2'-2k_6/3, \quad k_3= k_3'-k_6,
  \ee
  such that for $k_6\in 3\mathbb{Z}+g$, $g=0,1,2$,
\be 
\nu^{||}_g=-\frac{g}{3}\left(\begin{array}{c} 1 \\ 2 \\
                          0 \end{array}\right) \mod \mathbb{Z}.
\ee
The sum then evaluates to
\be
\begin{split}
&\sum_{k_j\in \mathbb{Z}}\kappa_3\, M_1(2\sqrt{2y/3}(k_6-9a))\, q^{\ul
  Q(\ul k)}e^{2\pi i\ul B(\ul k,\ul z)} =\\
& \times \sum_{k_6\in 3\mathbb{Z}+g \atop g=0,1,2}
                                                   q^{-2k_6^2/3}w^{12k_6}\,M_1(2\sqrt{2y/3}(k_6-9a))
                                           \,\Phi^+_{0,\nu^{||}_g}(\tau,
                                           u,v,\{d_1,d_2\}).
                                         \end{split}
                                         \ee
                                         As before, we relate
                                         $\Phi^+_{\nu^{||}_g}$ to an Appell
                                         function for $a$ as in Eq. (\ref{eq:-1a0}). Then with \be
\tilde \nu^{||}_0=0,\qquad 
\tilde \nu^{||}_1= \frac{1}{3}\left(\begin{array}{c} -1 \\ 1 \\
                          0 \end{array}\right), \qquad 
                      \tilde \nu^{||}_2= \frac{1}{3}\left(\begin{array}{c} 1 \\ 2 \\
                          0 \end{array}\right),
                      \ee
$\Phi^+_{0,\nu^{||}_g}$ agrees with $\Phi_{0,\tilde \nu^{||}_g}$.

\end{enumerate}

The results agree with the general structural formula (\ref{eq:StructForm}).

\subsection*{Terms with $L=2$}
Next we look at the summations corresponding to $\widehat K_2$ in
Eq. (\ref{eq:CRK123}). We write these as
\be
\label{eq:whatK2}
\begin{split}
\widehat K_2 &= \kappa'_1\, M_2\left(\frac{1}{\sqrt{2}}, 2\sqrt{y}(k_6-k_5/2-3a),\sqrt{2y}(k_5-12a) \right) \\ 
&\quad + \kappa'_2 \,M_2\left(\frac{1}{\sqrt{8}},\sqrt{3y}(k_4-k_6/3-6a),\frac{2\sqrt{2y}}{\sqrt{3}}(k_6-9a)\right) \\ 
&\quad + \kappa'_3 \,M_2\left( \frac{1}{\sqrt{2}}, 2\sqrt{y}(k_4-k_5/2-3a),\sqrt{2y}(k_5-12a) \right),
\end{split}
\ee
with
\begin{eqnarray}
\label{eq:kappaprimes}
\kappa'_1 &=& \tfrac{1}{2}(\sgn(-2k_1+k_2+4a)+\sgn(k_4-k_5/2-3a)),\\ \nn
\kappa'_2 &=& \tfrac{1}{2}(\sgn(k_1+k_2-k_3+4a)+\sgn(k_5-k_4/2-k_6/2-3a)),\\ 
\kappa'_3 &=& \tfrac{1}{2}(\sgn(2k_3-k_2+4a)+\sgn(k_6-k_5/2-3a)).\nn
\end{eqnarray} 
The $c_{i,j}$ read for these cases
\be
\begin{split}
  c_{3,1}=0,\qquad c_{3,2}=\frac{1}{2},\\
  c_{2,1}=\frac{1}{2},\qquad c_{2,3}=\frac{1}{2},\\
  c_{1,2}=\frac{1}{2},\qquad c_{1,3}=0.
\end{split} 
\ee 
\begin{enumerate}
\item For the term with $\kappa_1'$, we make the substitution
  $k_4=k'_4+k_5/2$ \eqref{eq:nu||}, and subsequently the substitution \eqref{eq:ktrafo}
\be
  k_1=k_1'-k_5/2,\quad k_2=k_2'-k_5, \quad
  k_3= k_3'-k_6.
  \ee
With the same arguments for $M_2$ as in Eq. \eqref{eq:whatK2}, the sum then becomes
  \be
  \begin{split}
&\sum_{k_j\in \mathbb{Z}} \kappa_1'\,M_2(\dots)\, q^{\ul Q(\ul k)/2}e^{2\pi i\ul B(\ul k,\ul
  z)} =\\
&\sum_{k_5\in 2\mathbb{Z}+g, g=0,1,\atop k_6\in \mathbb{Z}}
q^{-3k_5^2/4+k_5k_6-k_6^2}w^{9k_5+6k_6} M_2(\dots) \, \Phi^+_{0,\nu^{||}_g}(\tau,u,v,\{d_1\}), \label{k5tok6M2}
\end{split}
  \ee
  with
  \be
\nu^{||}_g=\frac{g}{2}\left(\begin{array}{c} 1 \\ 0 \\
                          0 \end{array}\right) \mod \mathbb{Z}.
  \ee
  
  The lattice $\Lambda_d(\{d_2^*,d_3^*\})$ is generated by $2d_3^*\pm
  d_2^*$, which matches the quadratic form for $k_5,k_6$ in this
  equation. Furthermore, 
  $|\Lambda_d^*(\{d_2^*,d_3^*\})/\Lambda_d(\{d_2^*,d_3^*\})|=8$, such
  that one finds for the
  number of glue vectors $\mathcal{N}_g=\sqrt{8\times 2/4}=2$. As
  before, we $\Phi^+_{0,\nu^{||}_g}$ matches with $\Phi_{0,\tilde \nu^{||}_g}$
  for $a$ as in Eq. (\ref{eq:-1a0}), with in this case
  \be
\tilde \nu^{||}_g=-\frac{g}{2}\left(\begin{array}{c} 1 \\ 0 \\
                          0 \end{array}\right).
  \ee
\item For the term with $\kappa_2'$, we substitute $k_5= k'_5+k_4/2+k_6/2$, and the
  substitution
\be
 k_1=k'_1+k_4/2-k_6/2, \quad
 k_2= k'_2-k_4/2-k_6/2, \quad k_3 = k'_3-k_6.
 \ee
The sum then evaluates to 
  \be
  \begin{split}
&\sum_{k_j\in \mathbb{Z}} \kappa_2'\,M_2(\dots)\, q^{\ul Q(\ul k)/2}e^{2\pi i\ul B(\ul k,\ul
  z)} =\\
&\sum_{k_4+k_6\in 2\mathbb{Z}+g \atop  g=0,1}
q^{-3k_4^2/4+k_4k_6-3k_6^2/4}w^{9k_4+9k_6} M_2(\dots) \, \Phi^+_{0,\nu^{||}_g}(\tau,u,v,\{d_2\}), \label{k4tok6M2}
\end{split}
  \ee
  with
  \be
\nu^{||}_g=\frac{g}{2}\left(\begin{array}{c} 1 \\ 1 \\
                          0 \end{array}\right) \mod \mathbb{Z}.
  \ee
   
  This matches with the general analysis. The lattice $\Lambda(\{d_1^*,d_3^*\})$ is generated by $d_1^*+d_3^*$
  and $4d_1^*$, with
  $|\Lambda^*(\{d_1^*,d_3^*\})/\Lambda(\{d_1^*,d_3^*\})|=8$. For
  the $\mathcal{N}_g=\sqrt{8\times 2/4}=2$ glue vectors we take $0$
  and $d_1$.

 For the same value of $a$ as in Eq. (\ref{eq:-1a0}), $\Phi^+_{0,\nu^{||}_g}$ matches with
 $\Phi_{0,\tilde \nu^{||}_g}$ with
 \be
\tilde \nu^{||}_g=\frac{g}{2}\left(\begin{array}{c} 1 \\ 1 \\
                          0 \end{array}\right).
  \ee
\item Finally for the term with $\kappa_3'$, we make the substitution
  $k_6= k'_6+k_5/2$, followed by the substitution
\be
k_1= k_1'+k_4-k_5, \quad k_2=k_2'-k_5, \quad
k_3= k'_3-k_5/2.
\ee
The sum then becomes
  \be
  \begin{split}
&\sum_{k_j\in \mathbb{Z}} \kappa_3'\,M_2(\dots)\, q^{\ul Q(\ul k)/2}e^{2\pi i\ul B(\ul k,\ul
  z)} =\\
&\sum_{ k_4\in \mathbb{Z}\atop k_5\in 2\mathbb{Z}+g, g=0,1}
q^{-k_4^2+k_4k_5-3k_5^2/4}w^{6k_4+9k_5} M_2(\dots) \, \Phi^+_{0,\nu^{||}_g}(\tau,u,v,\{d_3\}), \label{k4tok5M2}
\end{split}
  \ee
  with
  \be
\nu^{||}_g=\frac{g}{2}\left(\begin{array}{c} 0 \\ 0 \\
                          1 \end{array}\right) \mod \mathbb{Z}.
  \ee
As under Case 1., there are two glue vectors. As
  before, $\Phi^+_{0,\nu^{||}_g}$ matches with $\Phi_{0,\tilde \nu^{||}_g}$ for $a$ as in Eq. (\ref{eq:-1a0}), with in this case
  \be
\tilde \nu^{||}_g=\frac{g}{2}\left(\begin{array}{c} 0 \\ 0 \\
                          1 \end{array}\right).
  \ee

\end{enumerate}

\subsection*{Term with $L=3$}
  There is no geometric sum on the term with $L=3$ corresponding to
  $M_3$. We make the substitution \eqref{eq:ktrafo}
  \be
k_1=k_1'+k_4-k_5,\qquad k_2=k_2'-k_5,\qquad k_3=k_3'-k_6.
  \ee 
The sum over $k_j'$ equals $\Phi^+_{0,0}(\tau,-,v,\varnothing)$, which
is independent of $u$ and equals the
standard theta series for the $A_3$ root lattice. The sum becomes in this case 
\be
\begin{split}
&\sum_{k_j\in \mathbb{Z}} M_3(\dots)\, q^{\ul Q(\ul k)/2}e^{2\pi i\ul B(\ul k,\ul
  z)} =\\
&\sum_{ k_4,k_5,k_6 \in \mathbb{Z}}
q^{-k_4^2-k_5^2-k_6^2+k_4k_5+k_5k_6}w^{6(k_4+k_5+k_6)} M_3(\dots) \,
\Phi_{0,0}(\tau,-,v,\varnothing), \label{sumM3}
\end{split}
\ee

\subsection{$A_N$ Lattice} 
We discuss in this subsection the general case that $\Lambda$ is the $A_N$
root lattice. Rather than giving the full solution, we work out a few
aspects of Eq. \eqref{psi} with $r_i=1\;\forall \;i$, such that $\ell$
in Eq. \eqref{psi} equals $N$. The symbols $a$ and $b$ in this
subsection refer to the $a$ and $b$ as used in Eq. (\ref{psi}); in
particular $a$ is {\it not} related to the imaginairy part of the elliptic variable $z$
as in Eq. \eqref{eq:aula}.

\subsection*{$A_N$ Lattice: $\Psi_{(1,\dots,1),(0,0)}$}
We consider
first the case of $a=b=0$. We solve the constraint on the $b_i$ as
\begin{equation}
\sum_{i=1}^{N+1}b_i r_i=b=0,\implies \sum_{i=1}^N b_i=-b_{N+1}.
\end{equation}
We further change the summation variables from $b_i$ to $k_i$ as in Eq. \eqref{eq:tildeltol}.
This transforms the quadratic form for the $b_i$ (\ref{eq:tildeA}) to the standard $A_N$
quadratic form for the $k_i$ (\ref{eq:A_Nqform}). If we choose the
simple roots $\alpha_j$ as basis vectors, the vectors $d_j$ are related to
the $\alpha_j$ by the Weyl reflection in the hyperplane orthogonal to
$\alpha_1$. 

Comparing with Eq. \eqref{PhiIndef}, we evaluate the term with the elliptic variable $\ul B(\ul k, \ul z)$ as
\begin{eqnarray}
\ul B(\ul k,\ul z)=(2k_1-4k_2-2\sum_{j=3}^N k_j)z + 4z(\sum_{i=1}^N k_{N+i}),
\end{eqnarray}
where $z_j$ are the components of $\ul z=(z_1,z_2, \dots , z_N,
z_{N+1},\dots , z_{2N})^T$, which read
\begin{equation}
  \label{eq:zN}
\begin{split}
z_1 &= 2z(N-2),\\ 
z_j &= 2z j(N+1-j),\quad {\rm for}\; j=2,3,\dots ,N,\\ 
z_{N+j} &= -3z j(N+1-j), \quad {\rm for}\; j=1,2,3,\dots, N.
\end{split}
\end{equation}
This reduces to Eq. \eqref{eq:uzA3} for $A_3$. We define
$z_d=(z_{N+1},\dots, z_{2N})^T\in \Lambda \otimes \mathbb{C}$.

For generic $N$, let us consider the contribution of the term with
$L=N$ in the non-holomorphic part, ie the term involving $M_N$.
The change of basis to bring the quadratic form in block diagonal form
\eqref{eq:ChofB}, is the following substitution 
\begin{eqnarray}
k_1=k'_1+k_{N+1}-k_{N+2}, \quad k_j = k'_j-k_{N+j},\quad {\rm for}\; N\ge j>1.
\end{eqnarray}
If we write $\ul k=(k',k_d)\in \ul \Lambda$, with $k'=(k_1',\dots,
k_N')\in \Lambda$ and $k_d=(k_{N+1},\dots,k_{2N})^T\in \Lambda$
then $\ul Q(\ul k)=Q(k')-Q(k_d)$ with $Q(k)$ the $A_N$ quadratic
form. This substitution furthermore gives for $\ul B(\ul k,\ul z)$:
\begin{eqnarray}
\ul B(\ul k,\ul z)\rightarrow 6z\sum_{j=1}^N k_{N+j}+2zk'_1-4zk'_2-2z\sum_{j=3}^{N}k'_j,
\end{eqnarray}
such that
\be
\label{eq:MnSummand}
\begin{split} 
&q^{\ul Q(\ul k)/2}e^{2\pi i\ul B(\ul k,\ul z)} M_N(\{d_i^*\}, \Lambda ;\sqrt{2y}(k_d+{\rm Im}(z_d)/y)) \\
&\qquad = q^{-Q(k_d)/2} e^{2\pi i\,6z(\sum_{j=1}^N  k_{N+j})}
M_N( \{d_i^*\}, \Lambda ;\sqrt{2y}(k_d+{\rm Im}(z_d)/y) )\\
&\qquad \quad \times q^{Q(k')/2}e^{4\pi i z (k'_1-2k'_2-\sum_{j=3}^{N}k'_j)}.
\end{split}  
\ee
This reduces to Eq. \eqref{sumM3} for $N=3$. As discussed in App. \ref{AppWeylSym}, the Weyl group leaves the
quadratic form $Q$ invariant. As a result, sums of the summand
(\ref{eq:MnSummand}) over $k'$ and $k_d$ for different sets
$\{d_i^*\}$ can be equivalent, which reduces the complexity of the
non-holomorphic term. For the completion in Section \ref{sec:LatA3},
these correspond to the first and third term in Eqs (\ref{eq:kappas})
and (\ref{eq:kappaprimes}). 
 
\subsection*{$A_N$ Lattice: $\Psi_{(1,\dots,1),(a,b)}$}

Let us choose $a, b\in \mathbb{Z}$ not necessarily $0$ for $\Psi_{(1,\dots,1),(a,b)}$
\eqref{psi}. As before we choose all $r_i=1$ for $i=1,2,\dots ,N+1$ in
$SU(N+1)$, such that $r$ and $\ell$ in Eq. (\ref{psi}) equal
$N+1$. The effect of generic $a,b$ is that $k\in \Lambda + \mu-\nu$ for specific
$\mu$ and $\nu$, which we determine in the following. We have,
\begin{equation}
\sum_{i=1}^{N+1}b_i r_i=b,\implies \sum_{i=1}^N b_i=b-b_{N+1}.
\end{equation}
Using the Weyl reflection with respect to the root $\alpha_1$, we have
the following equations to relate the $b_i$ to $k_i$,
\begin{eqnarray}
b_1 &=& k_2-k_1+\frac{b}{N+1}, \\ \nn
b_2 &=& k_1+\frac{b}{N+1}, \\ \nn
b_j &=& k_j-k_{j-1}+\frac{b}{N+1}, \quad {\rm for}\; j=3,4,\dots ,N.
\end{eqnarray}
Since the $b_j\in\mathbb{Z}$, we have
\begin{eqnarray}\label{knu}
k_j &\in & \mathbb{Z}-\frac{jb}{N+1}.
\end{eqnarray}
Comparison with Eq. \eqref{phigen2} demonstrates that the elements of
$\mu=(\mu_1,\dots, \mu_N)^T$ and $\nu=(\nu_1,\dots, \nu_N)^T$ are given by
\be
\mu_j-\nu_j=-\frac{jb}{N+1} \mod \mathbb{Z}.
\ee
To determine $\mu$ and $\nu$ separately, we consider the terms in
Eq. (\ref{psi}) corresponding to $B(\nu,k)$ and $B(\nu,u)$. For
$B(\nu,k)$, this is
\be
\begin{split}
&\sum_{i=2}^{N+1}(b_{i-1}-b_i)\left\{\frac{a}{N+1}\sum_{j=i}^{N+1}1\right\}=\sum_{i=1}^N
(b_{i}-b_{i+1})\left\{ -\frac{ai}{N+1}\right\}\\
&\quad  = (k_2-2k_1) \left\{ -\frac{a}{N+1}\right\}+(k_1+k_2-k_3) \left\{
  -\frac{2a}{N+1}\right\}\\
&\qquad +\sum_{j=3}^{N-1}(-k_{j-1}+2k_j-k_{j+1})
\left\{ -\frac{ja}{N+1}\right\}\\
& \qquad +(2k_N-k_{N-1}) \left\{ -\frac{Na}{N+1}\right\}.
\end{split}
\ee
We find thus that $\nu$ is given
\be
\nu=\sum_{j=1}^N S_{\alpha_1}(\alpha_j) \left\{ -\frac{ja}{N+1}\right\},
\ee
or in components
\be
\nu=\left( \begin{array}{c} - \left\{ -\tfrac{a}{N+1}\right\} +\left\{
             -\tfrac{2a}{N+1}\right\}\\
             \left\{
             -\tfrac{2a}{N+1}\right\}  \\   \left\{
             -\tfrac{3a}{N+1}\right\} \\ \vdots \\ \left\{
             -\tfrac{Na}{N+1}\right\} \end{array}\right).
\ee
This agrees with the expression for $B(\nu,u)$ with $u=(z_1,\dots,z_N)^T$ \eqref{eq:zN},  
\begin{eqnarray}
4z \sum_{j=2 }^{N+1}\left\{\frac{a}{N+1}(N+2-j)\right\}&=&4z
                                                         \sum_{j=1}^{N}\left\{-\frac{ja}{N+1}\right\}.
\end{eqnarray}

\appendix

\section{The $A_N$ Root Lattice} 
\subsection{Roots of the $A_N$ Lattice}
\label{SULattices}
In this section we review some basic properties of the $A_N$ root lattice
of $SU(N+1)$ and its Weyl group. The root lattice contains $N$ simple
roots, $\alpha_1,\alpha_2,\cdots ,\alpha_N$, with unit norm with respect to the standard innerproduct $(\cdot)$ on
$\mathbb{R}^{N}$. The angle between
consecutive roots $\alpha_j, \alpha_{j\pm 1}$ is $\frac{2\pi}{3}$, and
the roots are orthogonal otherwise. The non-vanishing innerproducts are thus
\begin{eqnarray} 
  \label{RnInnerProduct}
\alpha_i^2=1,\quad \alpha_i\cdot \alpha_j=-\frac{1}{2}\delta_{i,j\pm 1}, \quad {\rm for}\; i\ne j.
\end{eqnarray}
Note this innerproduct is not even, and not integral. 
 
We furthermore introduce the Cartan matrix for $SU(N+1)$, which is the
even, integral quadratic form ${\bf A}_N$ with entries
\begin{eqnarray}
({\bf A}_N)_{ij} =\frac{2\alpha_i\cdot\alpha_j}{\alpha_j^2},
\end{eqnarray}
or
\be
\label{eq:A_Nqform}
{\bf A}_N = \left(\begin{array}{ccccc} 
    2 & -1 & 0 &  &  \\
    -1 & 2 & -1 & 0 & \\
   0 & -1 & \ddots & \ddots & \ddots\\
 & \ddots & \ddots &  \ddots &  1\\
                     & &    &   -1 & 2
    \end{array}\right), \quad \det({\bf A}_N)=N+1.
\ee

We let the quadratic form $Q(l)$ correspond to ${\bf A}_N$
(\ref{eq:A_Nqform}). Thus for $l=\sum_{i} l_i\alpha_i \in A_N$, $Q(l)$ reads 
 \be
 \label{defQ}
 Q(l)=2\sum_{i=1}^{N}l_i^2-2\sum_{i=1}^{N-1}l_il_{i+1}.
 \ee
 Furthermore, the bilinear form $B$ is defined as
 \be
B(l,l')=\frac{1}{4} (Q(l+l')-Q(l-l')),
 \ee
 such that $B(l,l)=Q(l)$, and $B(\alpha_i,\alpha_j)=({\bf A}_N)_{ij}$. 

The change of basis given by
 \be
 \label{eq:tildeltol}
 \begin{split}
\tilde l_1 &= l_2-l_1, \\ 
\tilde l_2 &= l_1, \\ 
\tilde l_j &= l_j-l_{j-1}, \quad {\rm for}\; j=3,4,\dots, N.
\end{split}
\ee 
changes the bilinear form ${\bf A}_N$ to
\be
\label{eq:tildeA}
\tilde {\bf A}_N = \left(\begin{array}{cccc} 
    2 & 1 & 1 & \dots \\
    1 & 2 & 1 & \dots\\
    \vdots & \ddots & & \\
     \dots & \dots & 1 & 2
    \end{array}\right),
\ee
with quadratic form $\tilde Q$:
 \be
\tilde Q(\tilde l)=2\sum_{i=1}^{N}\tilde l_i^2+\sum_{i\ne j}\tilde l_i\tilde l_{j}.
 \ee
\subsection{Weyl Group} 
\label{AppWeylSym}

The Weyl group of the root system of $SU(N+1)$ is the symmetric group $\mathcal{S}_{N+1}$, which
is generated by the reflections through the hyperplanes orthogonal to the $N$ simple roots
$\alpha_1,\dots, \alpha_{N}$. If $\rho\in A_N$ is a root, then under
the reflection $S_\alpha$ along the root
 $\alpha$ we have, 
\be
S_\alpha(\rho)= \rho-2\alpha \frac{ \rho\cdot \alpha  }{\alpha^2}.
\ee  This reflection preserves the inner product  $\rho_1\cdot \rho_2$,
\begin{eqnarray}\nn
  \rho_1\cdot \rho_2 &\rightarrow &  \rho_1\cdot \rho_2 - 2 \alpha\cdot \rho_2  \frac{
                          \rho_1\cdot\alpha  }{\alpha^2}-2 \rho_2\cdot\alpha
                          \frac{ \rho_1\cdot \alpha}{\alpha^2}+4 \rho_2\cdot
                          \alpha \frac{ \rho_1\cdot \alpha }{\alpha^2}=
                          \rho_1\cdot \rho_2. 
\end{eqnarray}
For the reflection $S_{\alpha_j}$ with respect to the simple root $\alpha_j$ at $j$-th position in the Dynkin diagram of $A_N$, we have the following transformation on the simple roots:
\begin{equation}\label{weylnj}
\begin{split}
&S_{\alpha_j}(\alpha_j)= -\alpha_j, \quad  S_{\alpha_j}(\alpha_{j-1})=\alpha_{j-1}+\alpha_j,\quad S_{\alpha_j}(\alpha_{j+1})=\alpha_{j+1}+\alpha_j,\\ 
&S_{\alpha_j}(\alpha_{j\pm k}) = \alpha_{j\pm k},\quad {\rm for}\;
k\ge 2.
\end{split}
\end{equation}
When expanded in terms of the basis $\{\alpha_i\}$, the image of a vector $k=\sum_i k_i\alpha_i
\in A_N$ under $S_{\alpha_j}$ is,
\begin{eqnarray}\label{bfk}
S_{\alpha_j}(k)=\sum_{i=1}^{j-1} k_i\alpha_i+\sum_{i=j+1}^N k_i\alpha_i -\alpha_{j}k_j+\alpha_{j}k_{j-1}+\alpha_{j}k_{j+1},
\end{eqnarray}
such that the components $k_j$ of $k$ are transformed as,
\begin{eqnarray}\label{kj}
k_{j}\rightarrow -k_j+k_{j-1}+k_{j+1}, \quad k_{i\ne j}=k_i.
\end{eqnarray} 
Clearly, the quadratic form $Q(k)$ (\ref{defQ}) remains invariant,
\be
\begin{split}
Q(k)/2 \rightarrow  & \sum_{i\ne j}k_i^2-\sum_{i\ne
  j,j-1}k_ik_{i+1}+(k_{j-1}+k_{j+1}-k_j)^2 \\
& -(k_{j-1}+k_{j+1})(k_{j-1}+k_{j+1}-k_j)\\ 
&= \sum_{i=1}^N k_{i}^2-\sum_{i=1}^{N-1}k_ik_{i+1}=Q(k)/2.
\end{split}
\ee

\section{Proof of the Orthogonality Relation}
\label{sec:ProofId}
In this Appendix we prove the orthogonality relation
\eqref{eq:orthorelation}. To this end, recall that by Eq. (\ref{eq:CsiPar}) the $c_{s,v}$ \eqref{eq:csv} are solutions to
$\ul B(C_{s_j}-\sum_{\ell=1}^L c_{s_j,v_\ell}C_{v_\ell},C_{v_k})=0$
for all $j=1,\dots, M-L$ and $k=1,\dots ,L$. Using
Eq. (\ref{eq:BCrCs}), this can be written in terms of the matrix $\mathbf{D}^{-1}$, 
\be 
\mathbf{D}^{-1}_{s_jv_k}=\sum_{\ell=1}^L c_{s_j,v_\ell}  \mathbf{D}^{-1}_{v_{\ell}v_k}.
\ee 
We multiply this equation by the matrix $\mathbf{D}_{s_is_j}$,
\be
\label{eq:OrthoRel}
\sum_{j=1}^{M-L} \mathbf{D}_{s_is_j}
\mathbf{D}^{-1}_{s_jv_k}=\sum_{j=1}^{M-L} \sum_{\ell=1}^L \mathbf{D}_{s_is_j} c_{s_j,v_\ell}  \mathbf{D}^{-1}_{v_{\ell}v_k}.
\ee 
Completing the sum over $s_j$ with $v_m$, $m=1,\dots,L$, we can write
the left hand side as
\be
\delta_{s_i,v_k}-\sum_{m=1}^L \mathbf{D}_{s_iv_m} \mathbf{D}^{-1}_{v_mv_k}.
\ee 
Since $s_i\neq v_k$ for all $i,k$, the $\delta$-function vanishes.
The identity (\ref{eq:OrthoRel}) therefore becomes
\be
\sum_{j=1}^{M-L} \sum_{\ell=1}^L (\mathbf{D}_{s_is_j} c_{s_j,v_\ell}+ \mathbf{D}_{s_iv_\ell} )\mathbf{D}^{-1}_{v_{\ell}v_k}=0.
\ee
Since ${\bf D}$ is positive definite, we can multiply this equation by
${\bf D}$ from the right. Substitution of ${\bf
  D}_{s,v}=B(d_{s},d_{v})$ demonstrates that this is equivalent to the desired identity
Eq. (\ref{eq:orthorelation}),
\be
B\left( d_{s_j},d_{v_\ell}+\sum_{k=1}^{M-L}c_{s_k,v_\ell}d_{s_k}\right)=0.
\ee

\providecommand{\href}[2]{#2}\begingroup\raggedright\endgroup


\end{document}